%% file: DwyerFinal10_28_07.tex
\documentclass{amsart}
\usepackage{amsmath,amsthm,amscd,amssymb}
\usepackage{pb-diagram}

\input amsmath-defs.tex
\newtheorem{thm}{Theorem}[section]
\newtheorem{lemma}[thm]{Lemma}
\newtheorem{prop}[thm]{Proposition}
\newtheorem{cor}[thm]{Corollary}
\theoremstyle{definition}

\newtheorem{defn}[thm]{Definition}

\newtheorem{rem}[thm]{Remark}

\numberwithin{equation}{section}

\newcommand{\hq}[1]{\ensuremath{H_#1(\underline{\ \ };\mathbb{Q})}}
\newcommand{\bz}{\mathbb{Z}}
\newcommand{\bq}{\mathbb{Q}}
\DeclareMathOperator{\im}{im}
\newcommand{\id}{\operatorname{id}}
\newcommand{\image}{\operatorname{image}}
\newcommand{\rank}{\operatorname{rank}}
\newcommand{\sss}{\scriptscriptstyle}
\newcommand{\2}[1]{\ensuremath{^{\sss (#1)}}}
\newcommand{\np}{\ensuremath{^{\sss (n+1)}}}
\newcommand{\nph}[1]{\ensuremath{#1^{\np}_{\sss H}}}
\newcommand{\nm}{\ensuremath{^{\sss (n-1)}}}
\newcommand{\gn}{\ensuremath{G^{\sss (n)}_{\sss H}}}
\newcommand{\gnp}{\ensuremath{G^{\sss (n+1)}_{\sss H}}}
\newcommand{\gone}{\ensuremath{G^{\sss (1)}_{\sss H}}}
\newcommand{\an}{\ensuremath{A^{\sss (n)}_{\sss H}}}
\newcommand{\bn}{\ensuremath{B^{\sss (n)}_{\sss H}}}
\newcommand{\cn}{\ensuremath{C^{\sss (n)}_{\sss H}}}

\title{Homology and Derived Series of Groups II: Dwyer's Theorem}
\author{Tim D. Cochran and Shelly Harvey$^\dag$}
\thanks{$^{\dag}$Both authors were partially supported by the National
Science Foundation; the second author was partially supported
by the Alfred P. Sloan Foundation}
\address{Rice University, Houston, Texas, 77005-1892}
\email{cochran@rice.edu, shelly@rice.edu}
\begin{document}
\begin{abstract}
We give new information about the relationship between the low-dimensional homology of a space and the derived series of its fundamental group. Applications are given to detecting when a set of elements of a group generates a subgroup ``large enough'' to map onto a non-abelian free solvable group, and to concordance and grope cobordism of classical links. We also greatly generalize several key homological results employed in recent work of Cochran-Orr-Teichner in the context of classical knot concordance.

In 1963 J. Stallings established a strong relationship between the low-dimensional homology of a
group and its lower central series quotients. In 1975 W. Dwyer extended Stallings' theorem by weakening the hypothesis on $H_2$. In 2003 the second author introduced a new characteristic series, $G\2n_{\sss H}$, associated to the derived series, called the
\emph{torsion-free derived series}. The authors previously established a precise analogue, for the torsion-free derived series, of Stallings' theorem. Here our main result is the analogue of Dwyer's theorem for the torsion-free derived series. We also prove a version of Dwyer's theorem for the \emph{rational} lower central series. We apply these to give new results on the Cochran-Orr-Teichner filtration of the classical link concordance group. 
\end{abstract}
\maketitle
\section{Introduction}\label{intro}

There are many situations in topology where the homology type of a space is fixed or is dependent only on coarse combinatorial data whereas the homotopy type, in particular the fundamental group, is a rich source of complexity. For example, for a knot $f:S^n\to S^{n+2}$ one sees by Alexander Duality that the exterior, $S^{n+2}-S^n$, is a homology circle, independent of the ``knotting'' of the embedding. Similarly, the homology groups of the exterior of an algebraic curve in $\mathbb{CP}(2)$ or $\mathbb{C}^2$ are determined merely by the intrinsic topology of the curve. Furthermore, in studying deformations of such embeddings, typically the homology groups of the exteriors do not vary, or are controlled by the combinatorics of the allowable singularities, whereas the fundamental group varies with few obvious constraints. Therefore to define interesting topological invariants of such embeddings, or of certain deformation classes of embeddings, it is vital to understand to what extent the homology of a space constrains its fundamental group. These issues are often profitably studied in purely group-theoretic terms, for if $X$ is a connected CW-complex then it is well-known that $H_1(X;\mathbb{Z})\cong H_1(\pi_1(X))$ and that the Hurewicz map induces an exact sequence
$$
\pi_2(X)\to H_2(X)\to H_2(\pi_1(X))\to 1.
$$
Thus group theoretic results quickly translate into results about spaces.

In 1963 John Stallings, in his landmark paper \cite{St}, established a strong relationship between the low-dimensional homology of a
group and its lower central series quotients. We review his results in abbreviated form. Recall that the $n^\supth$ term of the lower central series of
$G$, denoted $G_{n}$, is inductively defined by $G_1=G$,
$G_{n +1}=[G_{n},G]$. Stallings
also defines what we shall call the $\textbf{rational lower central
series}$, $G_{n}^r$, which differs only in that $G_{n
+1}^r$ consists of all those elements \emph{some finite power of
which} lies in $[G_{n}^r,G]$. 

\begin{thm}[{\cite[Theorem~3.4]{St}} Stallings' Integral Theorem]\label{Stallint}  Let $\phi:A\to B$ be a
homomorphism that induces an isomorphism on $H_1(-;\mathbb{Z})$ and an
epimorphism on $H_2(-;\mathbb{Z})$. Then, for each $n$, $\phi$ induces
an isomorphism $A/A_n\cong B/B_n$.
\end{thm}

\begin{thm}[{\cite[Theorem~7.3]{St}} Stallings' Rational Theorem]\label{Stallrat} Let $\phi:A\to B$ be a
homomorphism that induces an isomorphism on $H_1(-;\mathbb{Q})$ and an
epimorphism on $H_2(-;\mathbb{Q})$. Then, for each $n$, $\phi$
induces a monomorphism $A/A^r_{n}\subset B/B^r_{n}$, and induces isomorphisms
$(A^r_n/A^r_{n+1})\otimes \mathbb{Q}\cong (B^r_n/B^r_{n+1})\otimes
\mathbb{Q}$.
\end{thm}

Stallings' Rational Theorem has an elegant reformulation wherein the conclusion is replaced by the conclusion that $A$ and $B$ have the same \emph{Malcev completion}. This was made explicit in Bousfield ~\cite{B}. In 1975 William Dwyer extended Stallings' Integral theorem by weakening the hypothesis on $H_2$ and indeed found precise conditions for when a specific lower central series quotient was preserved (Theorem~\ref{DwyerTh1} below)~\cite{Dw}. For this purpose he defined an important subgroup of $H_2(A)$, denoted $\Phi_n(A)$, $n\geq 1$, as the kernel of $H_2(A)\to H_2(A/A_n)$. Dwyer's ``filtration'' of $H_2(A)$ has equivalent more geometric formulations in terms of surfaces (see below) and \emph{gropes}. These other formulations and the theorem below played a crucial role in Freedman and Teichner's work on $4$-manifold topology that strengthened the foundational results of Freedman-Quinn ~\cite{FT2}~\cite{FQ}~\cite{Kr1}~\cite{Kr2}~\cite{KT}.

\begin{thm}[{\cite[Theorem 1.1]{Dw}} Dwyer's Integral Theorem]\label{DwyerTh1}Let $\phi:A\to B$ be a
homomorphism that induces an isomorphism on $H_1(-;\mathbb{Z})$. Then for any positive integer $n$ the following are equivalent:
\begin{itemize}
\item $\phi$ induces
an isomorphism $A/A_{n+1}\cong B/B_{n+1}$.
\item $\phi$ induces an epimorphism $H_2(A;\mathbb{Z})/\Phi_n(A)\to H_2(B;\mathbb{Z})/\Phi_n(B)$. 
\end{itemize}
\end{thm}

For completeness, we herein prove an important result that is missing from the literature. Namely we use the Stallings-Dwyer methods to prove the following ``rational version of Dwyer's theorem'' for the (rational) lower central series. Here $\Phi_n^r(A)$, $n\geq 1$, is the kernel of $H_2(A;\mathbb{Z})\to H_2(A/A_n^r;\mathbb{Z})$. Indeed for any normal subgroup $N$ of $A$, we can define $\Phi_N(A)$ as the kernel of $H_2(A;\mathbb{Z})\to H_2(A/N;\mathbb{Z})$.
\newtheorem*{thm:DwyerTh2}{Theorem~\ref{DwyerTh2}}
\begin{thm:DwyerTh2}[Rational Dwyer's Theorem] Let $\phi:A\to B$ be a
homomorphism that induces an isomorphism on $H_1(-;\mathbb{Q})$. Then for any positive integer $n$ the following are equivalent:
\begin{itemize}
\item For each $1\leq k\leq n$, $\phi$ induces a monomorphism $A/A^r_{k+1}\subset B/B^r_{k+1}$, and isomorphisms
$H_*(A/A^r_{k+1};\mathbb{Q})\cong H_*(B/B^r_{k+1};\mathbb{Q})$; and an isomorphism $(A^r_k/A^r_{k+1})\otimes \mathbb{Q}\cong (B^r_k/B^r_{k+1})\otimes
\mathbb{Q}$.
\item $\phi$ induces an epimorphism $H_2(A;\mathbb{Q})/<\Phi_n^r(A)>\to H_2(B;\mathbb{Q})/<\Phi_n^r(B)>$. 
\end{itemize}
\end{thm:DwyerTh2}

Recall that the $n^\supth$
term of the derived series, $G^{\sss (n)}$, is defined by
$G\20=G$, $G^{\sss (n +1)}=[G^{\sss (n)},G^{\sss (n)}]$. Elementary examples show that the naive analogues of these theorems for the \emph{derived series} are false. For example If $A$ is the fundamental group of the exterior of a knot in $S^3$ then the abelianization map to $\mathbb{Z}$ is a homology equivalence, but $A^{(1)}/A^{(2)}$ is the Alexander module of the knot and is usually highly nontrivial. Therefore, until recently, there have been few if any links found between the homology of a group and its derived series. In 1974, Ralph Strebel had some partial success that is the starting point of our work ~\cite{Str}. However in ~\cite[Section 2]{Ha2} the second author introduced a new characteristic series,
$G\2n_{\sss H}$, associated to and containing the derived series, called the
\textbf{torsion-free derived series}. For free groups, the torsion-free derived series coincides with the derived series. Using this larger series, the authors established in ~\cite{CH1} a precise analogue of Stallings'
Rational Theorem (stated below). In order to state this result we review the definition of the torsion-free derived series. Observe that if a subgroup $G\2n_{\sss H}$ (normal in $G$) has been defined then
${G\2n_{\sss H}}/{[G\2n_{\sss H},G\2n_{\sss H}]}$ is not only an abelian group but also a right {\it module} over $\bz[G/G\2n_{\sss H}]$, where the action is induced from the
conjugation action of $G$ ($[x]g=[g^{-1}xg]$). Harvey was motivated (by the known failures of the quotients by the derived series to respect homological equivalences) to eliminate
torsion {\it in the module sense} from the successive quotients $G^{(n)}/G^{(n+1)}$. Specifically  set $G\20_{\sss H}=G$. Once
$G^{\sss (n)}_{\sss H}$ has been inductively defined, let $T_n$ be the subset
of ${\gn}/{[\gn,\gn]}$ consisting of the $\bz[G/G\2n_{\sss H}]$-torsion elements,
i.e. the elements $[x]$ for which there exists some non-zero
$\g\in\bz[G/G\2n_{\sss H}]$, such that $[x]\g=0$. (In fact, it was shown in \cite[Proposition 2.1]{Ha2} that $T_n$ is a \emph{submodule}). Consider the
epimorphism of groups:
$$
G^{\sss (n)}_{\sss H} \xrightarrow{\pi_n} \f{\gn}{[\gn,\gn]}
$$
and define $G^{\sss (n+1)}_{\sss H}$ to be the inverse image of $T_n$ under
$\pi_n$. Then $G^{\sss (n+1)}_{\sss H}$ is a normal subgroup
of $G\2n_{\sss H}$ that contains $[\gn,\gn]$. It follows inductively that
$G^{\sss (n+1)}_{\sss H}$ contains $G^{\sss (n+1)}$. Moreover,
since ${\gn}/{G\np_{\sss H}}$ is the quotient of the module ${\gn}/{[\gn,\gn]}$
by its
torsion submodule, it is a $\bz[G/G^{\sss (n)}_{\sss H}]$ torsion-free module
\cite[Lemma 3.4]{Ste}.
Hence the successive quotients of the
torsion-free derived subgroups are {\it torsion-free modules} over
the appropriate rings. We remark that if the successive quotients, $G^{(n)}/G^{(n+1)}$, of the derived series are torsion-free modules (as holds for a free group) then the torsion-free derived series coincides with the derived series.

\begin{thm}[{\cite[Theorem 4.1]{CH1}}]\label{harvey} Let $\phi:A\to B$ be a
homomorphism that induces a monomorphism on $H_1(- ;\mathbb{Q})$ and an
epimorphism on $H_2(- ;\mathbb{Q})$. Suppose also that $A$ is finitely-generated and $B$ is finitely related. Then, for each integer
$n$, $\phi$ induces a monomorphism
$A/A^{\sss (n)}_{\sss H} \subset B/B^{\sss (n)}_{\sss H}$. If, in addition, $\phi$
induces an isomorphism on $H_1(- ;\mathbb{Q})$ then
$\an/A\np_{\sss H} \to \bn/B\np_{\sss H}$ is a monomorphism between modules of the
same rank (over $\mathbb{Z}[A/\an]$ and $\mathbb{Z}[B/\bn]$, respectively). 
\end{thm}

The main purpose of the present paper is to prove an analogue of Dwyer's Theorem for the torsion-free derived series. In order to state the results, we first define the appropriate analogue of Dwyer's $\Phi_n$. The obvious analogue: the kernel of $H_2(A)\to H_2(A/A^{(n)})$ turns out to be the wrong one.

\begin{defn}\label{psin} Suppose $N$ is a normal subgroup of a group $A$. Let $\Phi^N(A)$ be the image of the inclusion-induced $H_2(N)\to H_2(A)$. Equivalently $\Phi^N(A)$ consists of those classes represented by maps of closed oriented surfaces $f:\Sigma\to K(A,1)$ such that $f_*(\pi_1(\Sigma))\subset N$. Such surfaces will be called \textbf{$N$-surfaces of $A$}. Specifically if $N=A^{(n)}$ we abbreviate $\Phi^N(A)$ by $\Phi^{(n)}(A)$, $n\geq 0$. Thus $\Phi^{(n)}(A)$ is the image of $H_2(A^{(n)})\to H_2(A)$, or, equivalently, the subgroup of $H_2(A)$ of elements that can be represented an oriented surface $f:\Sigma \to K(A,1)$ such that $f_*(\pi_1(\Sigma))\subset A^{(n)}$. If $N=A^{(n)}_H$ we abbreviate $\Phi^N(A)$ by $\Phi^{(n)}_H(A)$. Thus $\Phi^{(n)}_H(A)$ is the image of $H_2(A^{(n)}_H)\to H_2(A)$ or equivalently the subgroup of $H_2(A)$ of elements that can be represented an oriented surface $f:\Sigma \to K(A,1)$ such that $f_*(\pi_1(\Sigma))\subset A^{(n)}_H$. Note that since $A^{(n)}\subset A^{(n)}_H$, $\Phi^{(n)}(A)\subset \Phi^{(n)}_H(A)$. 
 
\end{defn}

From this definition, it may not be immediately apparent that this is a natural analogue, for the derived series, of Dwyer's $\Phi_n$, because even for $N=A_n$, $\Phi^{N}(A)$ is generally much smaller than $\Phi_{N}(A)$. To see in what sense $\Phi^{(n)}(A)$ is an analogue of $\Phi_{n}(A)$ first recall that Dwyer's subgroup has the following equivalent reformulation. For any space $X$, let $\Phi_n(X)$ be the subgroup of $H_2(X)$ consisting of those elements that can be represented by an oriented surface $f:\Sigma \to X$ such that, for some standard symplectic basis of curves $\{ a_i,b_i | 1\leq i \leq $genus$(\Sigma)\}$ of $\Sigma$, $f_*([a_i])\subset \pi_1(X)_n$. That is, \emph{one-half} of the symplectic basis of curves maps into $\pi_1(X)_n$. Observe that if $\pi_1(X)_n$ is killed then such homology classes become spherical. From this observation it is not difficult to see that Dwyer's $\Phi_n(A)$ is the same as $\Phi_n(K(A,1))$ in the sense of the reformulation (see ~\cite[Lemma 2.4]{FT2}). Then the correct analogue for the derived series is \emph{not} to merely replace $\pi_1(X)_n$ by $\pi_1(X)^{(n)}$ in the above definition, but rather to additionally require that a \emph{full} symplectic basis map into $\pi_1(X)^{(n)}$. Indeed $\Phi^{(n)}(A)$ as defined above is clearly  the subgroup of $H_2(A)$ consisting of elements represented by an oriented surface $f:\Sigma \to K(A,1)$ such that, for some standard symplectic basis of curves $\{ a_i,b_i | 1\leq i \leq$ genus$(\Sigma)\}$ of $\Sigma$, $f_*([a_i])\subset \pi_1(X)^{(n)}$ and $f_*([b_i])\subset \pi_1(X)^{(n)}$. 

For a low-dimensional topologist the following remark may be enlightening. Dwyer's $\Phi_n(X)$ is the subgroup consisting of those homology classes that can be \emph{represented by half-gropes of class n+1} whereas $\Phi^{(n)}(X)$ is the subgroup consisting of those classes that can be represented by \emph{symmetric gropes of height n}.

The following is then the main result of this paper.

\newtheorem*{thm:main}{Theorem~\ref{main}}
\begin{thm:main}[Main Theorem] Let $A$ be finitely-generated and $B$ finitely related.  Suppose $\phi:A\to B$ induces
a monomorphism on $H_1(- ;\mathbb{Q})$ and induces an epimorphism $\phi_*: H_2(A;\mathbb{Q})\to H_2(B;\mathbb{Q})/<\Phi^{(n)}_H(B)>$ (that is, the cokernel of $\phi_*: H_2(A;\mathbb{Q})\to H_2(B;\mathbb{Q})$ is spanned by $B^{(n)}_H$-surfaces). Then  $\phi$ induces a monomorphism $A/A_H^{(n+1)} \subset
B/B^{(n+1)}_H$. If, in addition, $\phi$ induces an isomorphism on
$H_1(- ;\mathbb{Q})$ then, in addition,
$\an/A\np_{\sss H}$ and $\bn/B\np_{\sss H}$ have the same rank (as modules over $\mathbb{Z}[A/\an]$ and
$\mathbb{Z}[B/\bn]$, respectively). 
\end{thm:main}

Since the torsion-free derived series of a free group is merely the ordinary derived series, we have the following application that \emph{makes no mention of the torsion-free derived series}.

\newtheorem*{thm:freegroup}{Corollary~\ref{cor:freegroup}}
\begin{thm:freegroup}Suppose $F$ is a free group, $B$ is a finitely-related group, $\phi:F\to B$ induces
a monomorphism on $H_1(- ;\mathbb{Q})$ and $H_2(B;\mathbb{Q})$ is spanned by $B^{(n)}$-surfaces. Then $\phi$ induces a monomorphism $F/F^{(n+1)} \subset
B/B^{(n+1)}$ (similarly for any $m\leq n+1$).
\end{thm:freegroup}

One of the algebraic applications of the work of Stallings and Dwyer is a criterion for when a set $\mathcal{A}= \{a_1,\dots,a_m\}$ of elements of a group $B$ generates a free subgroup $A$ of rank $m$. Indeed, if $B$ is itself a free group then it is a classical result that if $\mathcal{A}$ is linearly independent modulo $F_2$ then, for any $n$, it freely generates, in $F/F_n$, a free-nilpotent subgroup \cite[p.117,42.35,p.76,26.33]{Ne}. It follows that $A$ is free of rank $m$. Stallings improves on this result by weakening the hypothesis that $B$ be free to the hypothesis that $H_2(B;\mathbb{Q})=0$. (See also ~\cite[p.303-304]{Str}). The $H_2$ condition is still quite restrictive. Dwyer's work weakens this condition, replacing it by the condition that all Massey products of 1-dimensional classes vanish for $B$. More precisely, Dwyer shows that if the integral Massey products vanish up to and including order $n-1$ then $A$ is ``large enough'' that it maps onto the free nilpotent group of rank $m$ and nilpotency class $n$ ~\cite[Proposition 4.3]{Dw}. Our work can be applied to generalize these results to give a criterion for when $\mathcal{A}$ freely generates, in $B/B^{(n)}$, a free-\emph{solvable} subgroup, that is to say, $A$ is ``large enough'' to map onto the \emph{free solvable} group of rank $m$ and derived length $n$.

\newtheorem*{prop:freesolvable1}{Proposition~\ref{freesolvable1}}
\begin{prop:freesolvable1}Suppose that $B$ is a finitely-related group and $A$ is the subgroup generated by $\mathcal{A}=\{a_i | i\in \mathcal{I}\}\subset B$. Suppose $\mathcal{A}$ is linearly independent in $H_1(B;\mathbb{Q})$ and suppose that $H_2(B;\mathbb{Q})= <\Phi^{(n-1)}(B)>$. Then $A/A^{(n)}$ is the free solvable group of derived length $n$ on $\mathcal{A}$, that is, if $F$ is the free group on $\mathcal{A}$ then the map $F\to A$ induces an isomorphism $F/F^{(n)}\cong A/A^{(n)}$. In particular $A$ maps onto the free solvable group on $\mathcal{A}$ of derived length $n$ and hence is not nilpotent if $m>1$. Moreover $A/A^{(n)}$ embeds in $B/B^{(n)}$.
\end{prop:freesolvable1}

Stallings' theorems and Dwyer's extensions have also had many applications in topology and our results provide extensions of these. For example, if $L_0$ and $L_1$ are oriented, ordered, $m$-component links of circles in $S^3$, we say they are \emph{concordant} if there exist compact oriented annuli $\Sigma_i, 1\leq i\leq m$, properly and disjointly embedded in $S^3\times [0,1]$, restricting to $L_j$ on $S^3\times \{j\}$, $j=0,1$. By Alexander Duality, the inclusion maps $((S^3\times\{j\})-L_j)\to (S^3\times [0,1])- \coprod \Sigma_i)$ induce integral homology equivalences. Thus Stallings' theorem implies that the isomorphism type of each of the quotients $A/A_n$, where $A=\pi_1(S^3-L)$, is an invariant of the concordance type of a link. Using this, A. Casson showed that Milnor's invariants of links are concordance invariants ~\cite{Ca}. The invariance under concordance of other link invariants such as the rank of the Alexander module ($A^{(1)}/A^{(2)}$) and certain signature invariants were also established using Stallings' theorem ~\cite[p.52]{Hi}~\cite{Sm}~\cite{L2}~\cite{Fr1}. Our previous work on the torsion-free derived series ~\cite{Ha2}~\cite{CH1} has led to corresponding new ``higher-order'' concordance invariants arising from ``higher-order'' signatures and ranks ~\cite{Ha2}~\cite{Ha3}. The Dwyer's theorem on the lower-central series and our present work on the torsion-free-derived series allow one to show that these invariants are unchanged under weaker equivalence relations involving certain surfaces instead of the annuli that appear in the definition of link concordance. This allows for analogues of all of the above results. These other equivalence relations on knots and links have been much studied recently in many different contexts and are related to notions of \emph{gropes}~\cite{T}~\cite{KT}~\cite{ConT1}~\cite{ConT2}~\cite{COT}~\cite{COT2}~\cite{CT}~\cite{Con1}. In particular, as applications of our main theorem we show that a family of Cheeger-Gromov von Neumann $\rho$-invariants of links and 3-manifolds considered by Harvey in ~\cite{Ha2} are actually invariants of weaker equivalence relations involving gropes similar but \emph{more general} than those considered in ~\cite{COT}~\cite{CT} (generalizing ~\cite[Section 6]{Ha2}). In Section~\ref{applications} we extend the work of ~\cite[Section 6]{Ha2} on the Cochran-Orr-Teichner filtration $\mathcal{F}_{(n)}$ of the classical disk-link concordance group. For example our results allow for the following sharpening of Harvey's ~\cite[Theorem 6.8]{Ha2}. Definitions and details are given in Section~\ref{applications}.

\newtheorem*{thm:infgenerated}{Theorem~\ref{infgenerated}}
\begin{thm:infgenerated}In the category of m-component ordered oriented string links ($m>1$), each of the quotients $\mathcal{F}_{(n)}/\mathcal{F}^{\mathbb{Q}}_{(n.5)}$ contains a subgroup, consisting entirely of boundary links, whose abelianization has infinite $\mathbb{Q}$-rank.
\end{thm:infgenerated}

\section{The Main Theorem}\label{maintheorem}

In this section we recall and prove the main theorem.

\begin{thm}[Main Theorem]\label{main}  Let $A$ be finitely generated and $B$ finitely related.  Suppose $\phi:A\to B$ induces
a monomorphism on $H_1(- ;\mathbb{Q})$ and induces an epimorphism $\phi_*: H_2(A;\mathbb{Q})\to H_2(B;\mathbb{Q})/<\Phi^{(n)}_H(B)>$. Then $\phi$ induces a monomorphism $A/A_H^{(n+1)} \subset
B/B^{(n+1)}_H$ (similarly for any $m\leq n+1$). If $\phi$ induces an isomorphism on
$H_1(- ;\mathbb{Q})$ then $\phi$ induces a monomorphism 
$\an/A\np_{\sss H}\to\bn/B\np_{\sss H}$ between modules
of the same rank (over $\mathbb{Z}[A/\an]$ and
$\mathbb{Z}[B/\bn]$ respectively). 
\end{thm}

\begin{cor}\label{cor:freegroup} Suppose $F$ is a free group, $B$ is a finitely-related group, $\phi:F\to B$ induces
a monomorphism on $H_1(- ;\mathbb{Q})$ and $H_2(B;\mathbb{Q})$ is spanned by $B^{(n)}$-surfaces (or more generally $B^{(n)}_H$-surfaces). Then $\phi$ induces a monomorphism $F/F^{(n+1)} \subset
B/B^{(n+1)}$ (similarly for any $m\leq n+1$).
\end{cor}

Before proving Theorem~\ref{main}, we offer reasonable motivation for the hypothesis on $H_2$. Consider the epimorphism $\phi: A\to A/A^{(n+1)}\equiv B$. This homomorphism certainly induces isomorphisms $A/A^{(i)}\cong B/B^{(i)}$ for $0\leq i\leq n+1$. If an analog of Dwyer's theorem were to hold then $\phi$ would satisfy an appropriate condition on $H_2$. Therefore let us examine the cokernel of $H_2(A)\to H_2(A/A^{(n+1)})$ and allow it to suggest a reasonable condition on $H_2$. Thinking topologically, consider a normal generating set $\{\gamma_i\}$ for $A^{(n+1)}$. We may obtain an Eilenberg-Maclane space $K(A/A^{(n+1)},1)$ by adjoining to $K(A,1)$  $2$-cells $\Delta_i$ along circles $\partial(\Delta_i)$ representing the $\gamma_i$, and then adding cells of dimension $3$ and higher. Since $\gamma_i\in A^{(n+1)}$, $\gamma_i=\prod_j[\alpha_{ij},\beta_{ij}]$ where $\alpha_{ij},\beta_{ij}\in A^{(n)}$. Thus $\partial(\Delta_i)$ is the boundary of a map $f_i:\Sigma_i\to K(A,1)$ of a compact orientable surface $\Sigma_i$ with a standard symplectic basis $\{ a_{ij},b_{ij} | 1\leq j \leq \text{genus}(\Sigma_i)\}$ where $(f_i)_*([a_{ij}])=\alpha_{ij}$ and $(f_i)_*([b_{ij}])=\beta_{ij}$. The cokernel of $\phi_*:H_2(A)\to H_2(B)$ is generated by the set of closed surfaces $\{\Delta_i\cup f_i(\Sigma_i)\}$ which are $n$-surfaces of $B$ as defined above.

\begin{proof}[Proof of Theorem~\ref{main}] The proof of the first claim is by induction on
$n$. The case $n=0$ is clear since $A/A^{\sss (1)}_{\sss H}$ is merely
$H_1(A;\mathbb{Z})/\{\mathbb{Z}\text{-Torsion}\}$ and the hypothesis that $\phi$ induces
a monomorphism on $H_1(-;\mathbb{Q})$ implies that it also
induces a monomorphism on $H_1(-;\mathbb{Z})$ modulo torsion.
Now assume that the first claim holds for $n$, i.e. $\phi$ induces a monomorphism $A/\an \subset
B/\bn$. We will prove that it
holds for $n+1$, under the hypothesis that the cokernel of $\phi_*: H_2(A;\mathbb{Q})\to H_2(B;\mathbb{Q})$ is spanned by $B^{(n)}_H$-surfaces.

It follows from \cite[Proposition 2.2]{Ha2} that $\phi(A\np_{\sss H})\subset B\np_{\sss H}$. Hence the diagram below exists and is commutative. By the Five Lemma, it suffices to show that
$\phi$ induces a monomorphism $\an/A\np_{\sss H}\to\bn/B\np_{\sss H}$.
\begin{equation*}
\begin{CD}
1      @>>>    \f{\an}{A\np_{\sss H}}  @>>>   \f A{A\np_{\sss H}}  @>>>
\f A{\an}    @>>> 1\\
&&     @VV \phi V        @VV\phi_{n+1}V       @VV\phi_nV\\
1      @>>>    \f{\bn}{B\np_{\sss H}}  @>>>   \f B{B\np_{\sss H}}  @>>>
\f B{\bn}    @>>> 1\\
\end{CD}
\end{equation*}
The proof follows exactly the proof of ~\cite[Theorem 4.1]{CH1} until we reach our Proposition~\ref{prop:main}. However, we will present a slightly different approach, as suggested by the referee. The strategy is to show that, essentially by definition, the module $\f{\an}{A\np_{\sss H}}$ can be formulated in terms of the homology of $A$ with certain twisted coefficients, and to then to prove the key theorem showing that an integral homology equivalence guarantees a homology equivalence with twisted coefficients modulo torsion (Proposition~\ref{prop:main}).
For simplicity we abbreviate $A/\an$ by $A_n$ and $B/\bn$ by
$B_n$. 

Suppose that $\an/A\np_{\sss H}\to\bn/B\np_{\sss H}$ were \emph{not} injective. Then, by examining the diagram below, we see that there would exist an $a\in\an$ representing a \emph{non}-torsion class $[a]$ in ${\an}/{[\an,\an]}$ such that $\phi (a)$
represents a $\mathbb{Z} B_n$-torsion class, $[\phi(a)]$, in ${\bn}/{[\bn,\bn]}$.

\begin{equation*}
\begin{CD}
\f{\an}{[\an,\an]}          @>\pi_A>>   \f{\an}{A\np_{\sss H}}\\
@VV\phi V             @VV\phi V\\
\f{\bn}{[\bn,\bn]}          @>\pi_B>>   \f{\bn}{B\np_{\sss H}}\\
\end{CD}
\end{equation*}

Thus it suffices to show:
\begin{cor}\label{cor:homology}: Under the map $\phi_*:\f{\an}{[\an,\an]}\to \f{\bn}{[\bn,\bn]}$ if $[\phi(a)]$is a $\mathbb{Z}B_n$-torsion class then $[a]$ is a $\mathbb{Z}A_n$-torsion class.
\end{cor}

This Corollary has a homological interpretation. 

\begin{rem}\label{observation} $\an/[\an,\an]\cong H_1(\an;\mathbb{Z})\cong H_1(A;\mathbb{Z}[A/\an])$. The second equivalence is, for an algebraist, a consequence of Shapiro's Lemma ~\cite[Proposition 2.6,p.73]{BR}. For a topologist, $H_1(A;\mathbb{Z}[A/\an]
)$ may be thought of as the first homology with twisted coefficients of an aspherical space $K(A,1)$ where $\pi_1(K(A,1))\cong A$ and the coefficient system is induced by $\pi_1(K(A,1))\cong A\to A/\an$ \cite[p.335]{HS}. Then $H_1(K(A,1);\mathbb{Z}[A/\an])$ can be interpreted as the first homology module of the covering space of K(A,1) corresponding to the subgroup $\an$, which is $\an/[\an,\an]$ \cite[Theorems VI3.4 and 3.4*]{Wh}.
 \end{rem} 

In general the set of torsion elements (in the usual sense) of a module over an arbitrary ring is not a submodule; but over an Ore domain it is a submodule. Note that the solvable groups $A_n$ and $B_n$ are $\mathbb{Z}$-torsion free since the successive quotients of the torsion-free derived series are torsion free modules. Thus $\mathbb{Z} A_n$ and
$\mathbb{Z} B_n$ are right Ore domains (see \cite[Proposition 2.1]{Ha2}) and consequently they admit (and embed into) classical right
rings of quotients (which are division rings) $\SK(A_n)$ and $\SK(B_n)$, respectively \cite[p.591-592,p.611]{P}. Moreover our inductive hypothesis is that $\phi$ induces a
monomorphism $A_n\to B_n$ and hence a ring monomorphism
$\mathbb{Z} A_n\to\mathbb{Z} B_n$. Thus $\phi$ induces a monomorphism $\SK(A_n)\to\SK(B_n)$, which endows
$\SK(B_n)$ with a $\SK(A_n)-\SK(B_n)$ bimodule structure. Recall that any module over a division ring is free (a vector space). Thus any module $M$ over an Ore Domain $R$ has a well-defined rank
which is defined to be the rank of the vector space $M
\otimes_{R}\mathcal{K}(R)$
\cite[p.48]{Co}. Alternatively the rank can be defined to be the
maximal integer $m$ such that $M$ contains a submodule isomorphic
to $R^m$. We also recall that for any Ore domain $R$, its quotient field $\mathcal{K}(R)$ is a flat $R$-module ~\cite[Proposition 3.5]{Ste}.

In summary we have the following homological interpretation of the torsion-free derived series.

\begin{prop}\label{observation2}~\cite[Proposition 2.12]{CH1}
 \item 1) $\an/A^{\sss (n+1)}_{\sss H}$ is equal to $H_1(A;\mathbb{Z}[A/\an])$ modulo its $\mathbb{Z}[A/\an]$-torsion submodule.
 \item 2) $A^{\sss (n+1)}_{\sss H}$ is the kernel of the composition:
$$
A^{\sss (n)}_{\sss H} \xrightarrow{\pi_n} \f{\an}{[\an,\an]}=H_1(A;\mathbb{Z}[A/\an])\to H_1(A;\mathbb{Z}[A/\an])\otimes_{\mathbb{Z}[A/\an]}\SK (A/\an).
$$
\end{prop}
\begin{proof} Property $1)$ follows directly from Remark~\ref{observation} and the definition of $A^{\sss (n+1)}_{\sss H}$. For Property $2)$, note that tensoring with the quotient field $\SK (A/\an)$ kills precisely the $\mathbb{Z}[A/\an]$-torsion submodule \cite[Corollary.II.3.3]{Ste}.
\end{proof}

\begin{proof} [Proof of Corollary~\ref{cor:homology}] Consider the following commutative diagram 
\begin{equation*}
\begin{CD}
\f{\an}{[\an,\an]}@>\cong>>
H_1(A;\mathbb{Z} A_n)\\
@VV\id\ox\phi V @VV\id\ox\phi V\\
\f{\an}{[\an,\an]}\ox_{\mathbb{Z} A_n}\mathbb{Z} B_n @>\cong>>
H_1(A;\mathbb{Z} B_n)@>i_*>> H_1(A;\mathbb{Q} B_n)\\
@VV\phi\ox\id V  @VV\phi_* V @VV\phi_* V\\
\f{\bn}{[\bn,\bn]}@>\cong>>
H_1(B;\mathbb{Z} B_n)@>i_*>> H_1(B;\mathbb{Q} B_n)\\
\\
\end{CD}
\end{equation*}
where the horizontal isomorphisms follow from Remark~\ref{observation} and from the fact that since
$A_n\to B_n$ is a group monomorphism, $\mathbb{Z}B_n$ is a free, hence flat, $\mathbb{Z}A_n$ module.
We now consider the central vertical composition above $\psi=\phi_*\circ(\id\ox\phi)$. Let $\bar a=(\id\ox\phi)([a])$.
Suppose that $\psi([a])$ is $\mathbb{Z} B_n$-torsion. Then there is a non-zero $\delta\in \mathbb{Z} B_n$ such that $\phi_*(i_*(\bar a)\delta))=0$ in $H_1(B;\mathbb{Q} B_n)$. By Proposition~\ref{prop:main} below (proof postponed), applied with $B_0=\bn$ and $\G=B/\bn=B_n$, $i_*(\bar a)\delta$ is $\mathbb{Q} B_n$-torsion in $H_1(A;\mathbb{Q} B_n)$. 

\begin{prop}\label{prop:main} Suppose $\phi:A\to B$ induces a monomorphism $H_1(A;\bbq)\rightarrowtail H_1(B;\bbq)$ and $H_2(B;\bbq)$ is generated by the images of
$$
\phi_*: H_2(A;\bbq)\lra H_2(B;\bbq)\quad\text{and}\quad\text{incl}_*: H_2(B_0;\bbq)\lra H_2(B;\bbq),
$$
where $B_0$ is a normal subgroup of $B$ such that $B/B_0$ is PTFA. If $A$ is finitely-generated, $B$ is finitely-related and $\G=B/B_0$, the kernel of $H_1(A,\bbq\G)\to H_1(B,\bbq\G)$ is a $\bbq\G$-torsion module.
If $H_1(A,\bbq)\to H_1(B,\bbq)$ is surjective, then the cokernel of $H_1(A,\bbq\G)\lra H_1(B,\bbq\G)$ is a $\bbq\G$-torsion module.
\end{prop}
\noindent Hence $i_*(\bar a)$ is $\mathbb{Q} B_n$-torsion. So in fact (after multiplying by an integer) there is a non-zero element $\gamma'\in \mathbb{Z} B_n$ such that $i_*(\bar a)\gamma'=i_*(\bar a\gamma')=0$ in $H_1(A;\mathbb{Q} B_n)$. Since the kernel of $i_*$ is $\mathbb{Z}$-torsion, $\bar a\gamma'$ is $\mathbb{Z}$-torsion and so $\bar a$ is $\mathbb{Z} B_n$-torsion \cite[Corollary.II.3.3]{Ste}. Let $\la\in \mathbb{Z}B_n$ be a non-zero element that annihilates $\bar a$. By hypothesis, $A_n$ is (isomorphic to) a subgroup of $B_n$. This implies that $\bbz B_n$, viewed as a $\bbz A_n$-module, is a free module --- say with basis $\ST\subset B_n$ and $1\in\ST$ --- and that $\bbz A_n$ may be considered  a direct summand of $\bbz B_n$ and that $H_1(A;\mathbb{Z} A_n)$ can be viewed as a direct summand of $H_1(A;\mathbb{Z} B_n)$. Express $\la$ as a linear combination $\sum_{t\in\ST}\la_t\cd t$. For some $t$, $\lambda_t\neq 0$. Right multiply by $t^{-1}$ and use that every right multiple of $\la$ by an element of $B_n$ also annihilates $\bar a$ to conclude that $\bar a$ is annihilated by some $\sum_{t\in\ST}\la_t'\cd t$ wherein with $\la'_1$ is a non-zero element of $\mathbb{Z} A_n$. Thus
$$
0=[a]\cd\lambda_1'+ \sum_{t\neq 1}[a]\la_t'\cd t.
$$
Since each component must vanish, $0=[a]\cd\lambda_1'$, showing that $[a]$ is $\mathbb{Z} A_n$-torsion.
This concludes the proof of Corollary~\ref{cor:homology}.
\end{proof}

This finishes the inductive step of the proof of the first
part of Theorem~\ref{main}, modulo the proof of Proposition~\ref{prop:main}.

The second claim of Theorem~\ref{main} follows in an identical fashion as in ~\cite[Proof of Theorem 4.1]{CH1}. For assume that $\phi$ induces an isomorphism on
$H_1(-;\mathbb{Q})$. We must show that
$\an/A\np_{\sss H}\to\bn/B\np_{\sss H}$ is a monomorphism between modules
of the same rank. The fact that this is a monomorphism follows
from the first part of the theorem. Since $\an/A\np_{\sss H}$ and
$\an/[\an,\an]$ differ only by $\mathbb{Z}A_n$-torsion, they have
the same rank, $r_A$, as $\mathbb{Z}A_n$-modules. For the same
reason, $\an/A\np_{\sss H}\ox_{\mathbb{Z}A_n}\SK(B_n)$ and
$\an/[\an,\an]\ox_{\mathbb{Z}A_n}\SK(B_n)$ are isomorphic. By
\cite[Lemma 4.2]{CH1}, the former has $\SK(B_n)$-rank equal to $r_A$
and hence so does the latter, which we have identified with
$H_1(A;\SK(B_n))$. If $\phi$ induces an isomorphism on
$H_1(-;\mathbb{Q})$ then note that $B$ must be finitely
generated. Hence Proposition~\ref{prop:main} applies to show
that $H_1(A;\SK(B_n))\cong H_1(B;\SK(B_n))$. Thus the latter has
$\SK(B_n)$-rank equal to $r_A$. But by applying the same reasoning
as above, we see that it has $\SK(B_n)$-rank equal to $r_B$, the
$\mathbb{Z}B_n$-rank of $\bn/B\np_{\sss H}$.

Thus the entire proof of Theorem~\ref{main} is reduced to the proof of the Proposition~\ref{prop:main}. 

Instead of giving a ``topological proof'' of Proposition~\ref{prop:main} along the lines of ~\cite{CH1}, we will give a ``group homology'' proof suggested by the referee. Our original proof of Proposition~\ref{prop:main} can be viewed in the old versions of the paper on the xxx arxiv. Along the way it seems valuable to reprove, in this new manner, the corresponding crucial result of ~\cite{CH1}. We are grateful to the referee for suggesting this proof. 

\subsection*{2. Another Proof of Proposition 4.3 of [CH-I]} 

First we need several basic homological results. We could suppose more generally for the results below that $\Gamma$ is a locally indicable group, $R$ is a subring of the rationals, $\mathbb{Z}\subset R\subset \mathbb{Q}$, and that $R\Gamma$ is an Ore domain (hence admitting a classical skew field of quotients $\SK$). The most common situation under which these hypotheses are satisfied is when $\Gamma$ is a \emph{poly-(torsion-free-abelian) group} (PTFA group), that is when $\Gamma$ has a finite normal series $\Gamma_i$ wherein the successive quotients are torsion-free abelian groups ~\cite[p.305]{Str}\cite[Section 2]{COT}. Note that, for any group $B$, the quotient $B/B_H^{(n)}$ is a PTFA group, since the quotients $B_H^{(i)}/B_H^{(i+1)}$ are $\mathbb{Z}$-torsion-free.

Recall the following modest generalization, due to the authors, of a result of R. Strebel ~\cite[p.305]{Str}.

\begin{lemma}[{\cite[Lemma 4.4]{CH1}}]\label{strebel} Suppose $f:M\to N$ is a homomorphism
of right $R\G$-modules with $N$ projective. Let $\ov f=
f\otimes \id$ be the induced homomorphism of $R$-modules
$M\otimes_{R\G}R\to N\ox_{R\G}R$. Then
$\rank_{R\G}(\image f)\ge\rank_\mathbb{Q}(\image \ov f)$.
\end{lemma}

\begin{cor}\label{rank} Suppose $\SC_*$ is a projective right $R\G$ chain complex with $C_p$ finitely generated. Then
$$
\rank_{R\G} H_p(\SC_*)\le\rank_{\mathbb{Q}} H_p(\SC_*\otimes_{R\G}R).
$$
\end{cor}
\begin{rem} This Corollary is false if $C_p$ is not finitely generated.
\end{rem}
\begin{proof}[Proof of Corollary~\ref{rank}] Let $r_p= \rank_{R\G}C_p$. Since $C_p$ is finitely generated, $r_p$ is finite. Let $\{\ov\SC_*\}=\{\ov C_p,\ov\p_p\}=\{C_p\otimes_{R\G}R,\p_p\ox\id\}$. We claim that $r_p= \rank_{\mathbb{Q}}\ov C_p$. This is clear if $C_p$ is free. If $C_p$ is merely projective it requires a short argument that we leave to the reader. Now observe
\begin{align*}
\rank_{R\G}H_p(\SC_*)  &= \rank_{R\G}(\ker\p_p) - \rank_{R\G}(\image\p_{p+1})\\
&= r_p - \rank_{R\G}(\image\p_p) - \rank_{R\G}(\image\p_{p+1})\\
&\le r_p - \rank_\mathbb{Q}(\image\ov\p_p) - \rank_\mathbb{Q}(\image\ov\p_{p+1})\\
&= \rank_\mathbb{Q}(\ker\ov\p_p) - \rank_\mathbb{Q}(\image\ov\p_{p+1})\\
&= \rank_\mathbb{Q} H_p(\ov\SC_*),
\end{align*}
where the inequality follows from two applications of Lemma~\ref{strebel} above.
\end{proof}

\begin{prop}\label{prop:oldCH} (~\cite[Proposition 4.3]{CH1}) Suppose $\phi:A\to B$ induces a monomorphism $H_1(A,\bbq)\rightarrowtail H_1(B,\bbq)$ and an epimorphism $H_2(A,\bbq)\to H_2(B,\bbq)$. If $A$ is finitely-generated, $B$ is finitely-related and $\G=B/B_0$ is PTFA, the kernel
of $H_1(A,\bbq\G)\lra H_1(B,\bbq\G)$ is a $\bbq\G$-torsion module. If, in addition, $H_1(A,\bbq)\rightarrowtail H_1(B,\bbq)$ is surjective, then the cokernel of $H_1(A,\bbq\G)\lra H_1(B,\bbq\G)$ is a $\bbq\G$-torsion module. (Since $\mathcal{K}\G$ is flat, the conclusions are the same as saying that $H_1(A,\mathcal{K}\G)\to H_1(B,\mathcal{K}\G)$ is a monomorphism (respectively an isomorphism)).
\end{prop}

\begin{proof}[Proof of Proposition~\ref{prop:oldCH}]  Choose a $\bbz A$-free resolution $X_*\thra\bbz$ and a $\bbz B$-free resolution $Y_*\thra\bbz$ with $X_1$ and $X_0$ finitely-generated over $\bbz A$ and $Y_2$ finitely-generated over $\bbz B$. Use $\phi:A\to B$ to turn right $\bbz B$-modules into right $\bbz A$-modules and lift the identity $\mathbf{1}_{\mathbb{Z}}$ to a chain map $\phi_*:X_*\to Y_*$ of $\bbz A$-modules. Next pass from the $\bbz A$ complex $X_*$ to the $\bbz B$-complex $\wt X_*=X_*\ox_{\bbz A}\bbz B$; then form the mapping cone $Z_*$ ~\cite[p.6]{BR}. The chain group in dimension $i$ of this complex is the free $\bbz B$-module $Z_i= \wt X_{i-1}\op Y_i$. For every left $\bbz B$-module $M$ the homology groups of the mapping cone $Z_*$ fit into the exact sequence ~\cite[Proposition 0.6]{BR}
\begin{equation}\label{eq:sequence}
\begin{split}
\dots  &\lra H_2(A,M)\lra H_2(B,M)\lra H_2(Z_*\ox_{\bbz B} M)\lra\\
&\lra H_1(A,M)\lra H_1(B,M)\lra H_1(Z_*\ox_{\bbz B} M)\lra
\end{split}
\end{equation}
The hypotheses imply that $H_*(Z_*\ox_{\bbz B} \mathbb{Q})$ vanishes in dimension 2, respectively in dimensions 2 and 1. Applying Corollary~\ref{rank} to the complex $Z_*\ox_{\bbz B} \bbq\G$ in dimension 2, respectively in dimensions 2 and 1, we conclude that $H_*(Z_*\ox_{\bbz B} \bbq\G)$ has $\bbq\G$-rank zero, which is to say, is $\bbq\G$-torsion in these same dimensions. To apply Corollary~\ref{rank} we need to know that $Y_1$ and hence $Z_1$ is finitely-generated. But note that if $\phi_*$ is an \emph{isomorphism} on $H_1(-;\mathbb{Q})$, then since $A$ is finitely generated, $H_1(B;\mathbb{Q})$ is finite-dimensional. Since $B$ is finitely related this can only happen if $B$ is finitely presented. Therefore in this case we can assume both $Y_1$ and $Y_2$ are finitely generated. Then the sequence ~\ref{eq:sequence} implies Proposition~\ref{prop:oldCH}. 

We note in passing that this same type of analysis of ~\ref{eq:sequence} shows that the cokernel part of Proposition~\ref{prop:oldCH} is true \emph{without any} hypothesis on $H_2$, as long as either $\phi_*$ is an isomorphism on $H_1(-;\mathbb{Q})$, or $\phi_*$ is an epimorphism on $H_1(-;\mathbb{Q})$ and $B$ is finitely presented (see ~\cite[Proposition 2.10]{COT}).
\end{proof}
 
\subsection*{3. Proof of Proposition~\ref{prop:main}} 

In the remainder of this proof all tensor products are over $\bbz B$ unless specified otherwise. The strategy of the proof is similar to that of Proposition~\ref{prop:oldCH}, but there is an new difficulty: the hypotheses no longer imply that $H_2(Z_*\ox\bbq)=0$ in the exact sequence ~\ref{eq:sequence} and so the desired conclusion does not follow directly from Corollary~\ref{rank}. However, if $H_1(A,\bbq)\rightarrowtail H_1(B,\bbq)$ is surjective, then the cokernel of $H_1(A,\bbq\G)\lra H_1(B,\bbq\G)$ is a $\bbq\G$-torsion module, since this part of the proof of Proposition~\ref{prop:oldCH} was not dependent on the $H_2$ condition. Thus the cokernel part of Proposition~\ref{prop:main} holds.

Here is a sketch of the proof. The kernel of $H_{1}(A;\bbq\Gamma)\rightarrow H_{1}(B;\bbq\Gamma)$ is a torsion module precisely when the cokernel of $H_{2}(B;\bbq\Gamma) \rightarrow H_{2}(B,A;\bbq\Gamma)$
is a torsion module, so we shall establish the latter. We observe that $H_{2}(B_{0};\bbq)$ can be canonically identified with $H_{2}(B;\bbq[B/B_{0}])$. With this in mind, our hypotheses ensure that $H_{2}(B;\bbq)\rightarrow H_{2}(B,A;\bbq)$
is surjective and moreover that any class in $H_{2}(B,A;\bbq)$ lies in the image of a class in $H_{2}(B;Q)$ that ``lifts'' to $H_{2}(B_{0};\bbq)=H_{2}(B;\bbq\Gamma)$.  Then we show that these lifted
classes span a submodule of full rank in $H_{2}(B,A;\bbq\Gamma)$, i.e. one whose $\bbq\Gamma$-rank equals that of $H_{2}(B,A;\bbq\Gamma)$. Corollary~\ref{rank} is used to ensure that the rank of $H_{2}(B,A;\bbq\Gamma)$
is controlled by that of $H_{2}(B,A;\bbq)$, and Lemma~\ref{strebel} is used to show that the lifted classes span modules of the same $\bbq\Gamma$-rank as the $\bbq$-rank of the span of their images in $H_{2}(B,A;\bbq)$. The following makes these ideas precise.

In other words our goal is to prove that the cokernel of $H_{2}(Y_{\ast} \otimes \bbq\Gamma)\rightarrow H_{2}(Z_{\ast} \otimes \bbq\Gamma)$ is a torsion module. Consider the inclusion $i_2:Y_2\to Z_2$ and the boundary map $\partial_3:Z_3\to Z_2$ and use them to define two homomorphisms
$$
(\partial_3\ox \mathbf{1}_{\bbq\G})\op (i_2\ox \mathbf{1}_{\bbq\G}): (Z_3\ox \bbq\G)\op(Y_2\ox\bbq\G)\lra Z_2\ox \bbq\G,
$$
$$
(\partial_3\ox \mathbf{1}_{\bbq})\op (i_2\ox \mathbf{1}_{\bbq}): (Z_3\ox \bbq)\op(Y_2\ox\bbq)\lra Z_2\ox \bbq.
$$
By restricting in domain and range we arrive at
\begin{eqnarray}\label{eq:ggbar}
g : (Z_3\ox\bbq\G) \op \ker(\partial_2^Y\ox \mathbf{1}_{\bbq\G})\lra \ker(\partial_2^Z\ox \mathbf{1}_{\bbq\G}),\\
\bar{g} : Z_3\ox\bbq \op \ker(\partial_2^Y\ox \mathbf{1}_\bbq)\lra \ker(\partial_2^Z\ox \mathbf{1}_\bbq).
\end{eqnarray}
By hypothesis $H_1(X_*\ox\bbq)\rightarrowtail H_1(Y_*\ox\bbq)$ is injective so $H_2(Y_*\ox\bbq)\to H_2(Z_*\ox\bbq)$ is surjective.  It follows that that the map $\bar{g}$ is surjective. 
Since $Z_2$ is finitely generated, $\ker(\partial_2^Z\ox \mathbf{1}_\bbq)$ is a finite-dimensional vector space. We can choose a basis, $\{\bar u'_1,\dots,\bar u'_r,\bar v'_1,\dots,\bar v'_s\}$, for this vector space where $\{\bar u'_i\}$ is a basis for the image of $\partial_3\ox\mathbf{1}_\bbq$ (say $\bar u'_i=(\partial_3\ox\mathbf{1}_\bbq)(\bar u_i)$ for $\bar u_i\in Z_3\ox\bbq$); and $\bar v'_i= (\imath_2\ox \mathbf{1}_{\bbq})(\bar v_i)$ for some $\bar v_i\in \ker(\partial_2^Y\ox \mathbf{1}_\bbq)$. Moreover we can choose the $\bar v_i$ more precisely. For, by the hypothesis on $H_2(B;\mathbb{Q})$, $\ker(\partial_2^Y\ox \mathbf{1}_\bbq)$ is generated by three types of classes:
\begin{itemize}
\item [A.] the image of $~\partial_3^Y\ox \mathbf{1}_\bbq:Y_3\ox\bbq\to Y_2\ox\bbq$,
\item [B.] the image of (cycles) $~\phi_*:X_2\ox\bbq\to Y_2\ox\bbq$, and
\item [C.] the image (under inclusion) of the kernel of
$$
\partial_2\ox \mathbf{1}_\bbq: Y_2 \ox_{\mathbb{Z}B_0}\bbq\lra Y_1\ox_{\mathbb{Z}B_0}\bbq.
$$
\end{itemize}
For type $C$ we have used the fact that $Y_*$ can also be viewed as a free $\mathbb{Z}B_0$-resolution to compute $H_2(B_0;\bbq)$. We claim that we can choose our $\bar v_i$ to be of type $C$.  For elements of type $B$ vanish under $\imath_2\ox \mathbf{1}_{\bbq}$; and since $\imath_*\ox \mathbf{1}_\bbq$ is a chain map,
$$
(\imath_2\ox \mathbf{1}_{\bbq})(\partial_3^Y\ox \mathbf{1}_\bbq)=(\partial_3^Z\ox \mathbf{1}_\bbq)(\imath_3\ox \mathbf{1}_{\bbq}),
$$
so elements of type $A$ are carried under $\imath_2\ox \mathbf{1}_{\bbq}$ into the subspace spanned by $\{\bar u'_i\}$.

Now we lift the set $\{\bar u_1,\dots,\bar u_r,\bar v_1,\dots,\bar v_s\}$ to a set $\{u_1,\dots,u_r, v_1,\dots,v_s\}$ in the domain of $g$ as follows. The elements $\bar u_i$ can be lifted to elements $u_i\in Z_3\ox\bbq\G$. Since $Y_*\ox_{\mathbb{Z}B_0}\bbq$ is canonically isomorphic to $Y_*\ox_{\mathbb{Z}B}\bbq\G$, the elements $\bar v_i$ can likewise be lifted to elements $v_i$ in the kernel of

$$
\partial_2^Y\ox \mathbf{1}_{\bbq\G}:Y_2\ox\bbq\G\to Y_1\ox\bbq\G.
$$

These sets induce maps $f$ and $\bar f$ of free modules into the domains of $g$ and $\bar g$ respectively as shown in the diagram.
$$
\begin{diagram}\dgARROWLENGTH=1.2em
\node{\<u_i,v_i\>} \arrow{e,t}{\equiv}
\arrow{s,r}{\pi} \node{(\bbq\G)^{r+s}} \arrow{e,t}{f}
\arrow{s,r}{\pi}\node{\text{domain}(g)}
\arrow{e,t}{g} \node{\ker(\partial_2^Z\ox \mathbf{1}_{\bbq\G})}\arrow{e,t}{\subset} \node{Z_2\ox \bbq\G}\arrow{s,r}{\pi}\\
\node{\<\bar u_i,\bar v_i\>} \arrow{e,t}{\equiv}\node{(\bbq)^{r+s}} \arrow{e,t}{\bar f}\node{\text{domain}(\bar g)}\arrow{e,t}{\bar g} \node{\ker(\partial_2^Z\ox \mathbf{1}_\bbq)} \arrow{e,t}{\subset} \node{Z_2\ox \bbq} 
\end{diagram}
$$
Since $\{\bar u'_1,\dots,\bar u'_r,\bar v'_1,\dots,\bar v'_s\}$ is a basis for $\ker(\partial_2^Z\ox \mathbf{1}_\bbq)$, the bottom composition is a monomorphism. By Lemma~\ref{strebel} (actually Strebel's original result suffices here ~\cite[p.305]{Str}), the top composition is also a monomorphism. Hence $g \circ f$ is a monomorphism. Thus the images of $\{u_i,v_i\}$ generate a rank (r+s) free $\bbq\G$-submodule $F$ of $\ker(\partial_2^Z\ox \mathbf{1}_{\bbq\G})$,
so the latter has rank at least $r+s$. On the other hand,
$$
\rank_{\bbq\G}\ker(\partial_2^Z\ox \mathbf{1}_{\bbq\G})=\rank_{\bbq\G}(Z_2\ox \bbq\G)-~\rank_{\bbq\G}\im(\partial_2^Z\ox \mathbf{1}_{\bbq\G})
$$
which, by Lemma~\ref{strebel} is at most
$$
\rank_{\bbq\G}(Z_2\ox \bbq\G)-~\rank_{\bbq}\im(\partial_2^Z\ox \mathbf{1}_{\bbq})
$$
$$
=\rank_{\bbq}(Z_2\ox \bbq)-~\rank_{\bbq}\im(\partial_2^Z\ox \mathbf{1}_{\bbq})
$$
$$
=\rank_{\bbq}\ker(\partial_2^Z\ox \mathbf{1}_{\bbq})=~r+s.
$$
Thus the $\bbq\G$-rank of $\ker(\partial_2^Z\ox \mathbf{1}_{\bbq\G})$ is exactly $r+s$. Since the first $r$ elements of the basis of $F$ are in the image of $Z_3\ox \bbq\G\to Z_2\ox\bbq\G$ (i.e. are boundaries) and the remaining ones are images of classes in the kernel of $Y_2\ox \bbq\G\to Y_1\ox\bbq\G$, the image of $H_2(Y_*\ox\bbq\G)\to H_2(Z_*\ox \bbq\G)$ is a submodule of full rank. Thus the cokernel of this map is a torsion module, which was what was to be proved. Therefore the kernel of $H_1(A,\bbq\G)\to H_1(B,\bbq\G)$ is a torsion module.

 This completes the proof of Proposition~\ref{prop:main} and hence the proof of Theorem~\ref{main}.
  
\end{proof} 

The following slightly more general version of Proposition~\ref{prop:main} is useful in the applications. It differs only in that the coefficient system $\psi:B\to\G$ is not assumed to be surjective.

\begin{prop}\label{2-connected} Suppose $\phi:A\to B$ induces a monomorphism (respectively an isomorphism) on $\hq1$ with $A$ finitely-generated and $B$ finitely-related. Consider the coefficient system $\psi:B\to\G$ where $\G$ is PTFA. Let $B_0=\ker\psi$. Suppose $H_2(B;\bbq)$ is spanned by $\phi_*(H_2(A;\bbq))$ together with a collection of $B_0$-surfaces. Then $\phi$ induces a monomorphism (respectively, an isomorphism)
$$
\phi_*: H_1(A;\SK\G)\lra H_1(B;\SK\G).
$$
\end{prop}

\begin{proof}[Proof of Proposition~\ref{2-connected}] Let $\G'=\text{image}(\phi)\cong B/B_0$. Then $\G'\subset \G$ is a subgroup of a PTFA group, hence is a PTFA group. There is an induced monomorphism
$i:\mathcal{K}\G'\to \mathcal{K}\G$ so $\mathcal{K}\G$ is a free, hence flat $\mathcal{K}\G'$-module (indeed every $\mathcal{K}\G'$-module is free). Thus we have the commutative diagram
\begin{equation*}
\begin{CD}
H_1(A;\mathcal{K}\G')          @>i_*>>   H_1(A;\mathcal{K}\G) @>\cong>> H_1(A;\mathcal{K}\G')\ox_{\mathcal{K}\G'}\mathcal{K}\G\\
@VV\phi'_* V             @VV\phi_* V  @VV\phi'_*\ox \mathbf{1} V\\
H_1(B;\mathcal{K}\G')          @>i_*>>   H_1(B;\mathcal{K}\G) @>\cong>> H_1(B;\mathcal{K}\G')\ox_{\mathcal{K}\G'}\mathcal{K}\G\\
\end{CD}
\end{equation*}
By Proposition~\ref{prop:main}, $\phi_*'$ is injective (respectively an isomorphism). Since $\mathcal{K}\G$ is flat, $\phi'_*\ox \mathbf{1}$ is injective (respectively an isomorphism). 
\end{proof}

The following variation of Corollary~\ref{rank} is often useful in applications.  The seemingly tiny technical results, Corollary~\ref{rank} and Corollary~\ref{cor:independence}, actually greatly generalize the key homological results employed in recent work of Cochran-Orr-Teichner establishing new techniques in knot theory. Their key homological results \cite[Proposition 4.3]{COT} were proven only in the context where the chain complexes were cellular chain complexes for a $4$-dimensional manifold, so it is somewhat surprising that in fact (as shown by Corollary~\ref{rank} and Corollary~\ref{cor:independence}) those hypotheses turn out to be quite superfluous!

\begin{cor}\label{cor:independence} Suppose $\SC_*$ is a projective chain complex of (right) $R\G$-modules with $C_p$ finitely generated. Let $\ov\SC_*=\SC_*\otimes_{R\G}R$ and $\pi_*:\SC_*\to\ov\SC_*$ be the obvious chain map sending $x$ to $x\otimes 1$. Suppose $\{x_s\mid s\in S\}$ is a set of $p$-cycles of $\SC_*$. Let $\ov x_s=\pi(x_s)$ and suppose the $\mathbb{Q}$- rank of $H_p(\ov\SC_*)/\<[\ov x_s]\mid s\in S\>$ is at most $k$. Then the $\SK$-rank of $H_p(\SC_*)/\<[x_s]\mid s\in S\>$ is at most $k$.
\end{cor}

\begin{proof}[Proof of Corollary~\ref{cor:independence}] We shall define a projective chain complex $\SD_*=\{D_p,d_p\}$ such that $H_p(\SD_*)\cong H_p(\SC_*)/\<[x_s]\mid s\in S\>$, and then apply Corollary~\ref{rank} to this chain complex. Set $D_{p+1}=(\oplus_{s\in S}R\G)\bigoplus C_{p+1}$ and otherwise set $D_i=C_i$. Let $d_{p+1}:D_{p+1}\to D_p$ be defined by $d_{p+1}(e_s,y)=x_s+\p_{p+1}(y)$ where $\{e_s\}$ is a basis of $(R\G)^s$ and $y\in C_{p+1}$, and $d_{p+2}:D_{p+2}\to D_{p+1}$  by $d_{p+2}(z)=(0,\p_{p+2}(z))$ for $z\in D_{p+2}=C_{p+2}$. Then the $p$-cycles of $\SD_*$ are the same as those of $\SC_*$ while the group of $p$-boundaries is larger (includes $\{x_s\}$). Hence $H_p(\SD_*)\cong H_p(\SC_*)/\<[x_s]\mid s\in S\>$ as claimed.

It now suffices to show that the $\SK$-rank of $H_p(\SD_*)$ is at most $k$. Note that $D_p=C_p$ is finitely generated. Just as above we can create a chain complex $\ov\SD_*$ which agrees with $\ov\SC_*$ except in dimension $p+1$ where $\ov D_{p+1}=(\oplus_{s\in S}R)\bigoplus\ov\SC_*$ with $\bar d_{p+1}:\ov D_{p+1}\to\ov D_p$ given by $d_{p+1}(\bar e_s,\bar y)=\bar x_s+\ov\p_{p+1}(\bar y)$ for $\{\bar e_s\}$ a basis of $R^s$ and $\bar y\in\ov C_{p+1}$. Then just as above:
$$
H_p(\ov\SD_*)\cong H_p(\ov\SC_*)/\<[x_s] ~| s\in S\>,
$$
which has $\mathbb{Q}$-rank at most $k$ by assumption. Moreover the chain map $\pi:\SC_*\to\ov\SC_*$ extends to $\tl\pi:\SD_*\to\ov\SD_*$ by setting $\tl\pi(e_s,y)=(\bar e_s,\pi(y))$, and one sees that $\ov\SD_*=\SD_*\otimes_{R\G}R$. It follows immediately from Corollary~\ref{rank} that the $\SK$-rank of $H_p(\SD_*)$ is at most $k$.
\end{proof}

\section{A Rational Dwyer's Theorem for the Lower Central Series}\label{Dwyerrational}

For completeness, we prove the missing ``rational'' version of Dwyer's theorem, that is the correct generalization of Stallings' Rational Theorem. This follows from the Stallings'-Dwyer techniques but was not stated by Dwyer, nor by Bousfield ~\cite{B}.

\begin{thm}[Rational Dwyer's Theorem]\label{DwyerTh2} Let $\phi:A\to B$ be a
homomorphism that induces an isomorphism on $H_1(-;\mathbb{Q})$. Then for any positive integer $n$ the following are equivalent:
\begin{itemize}
\item [1)] For each $1\leq k\leq n$, $\phi$ induces a monomorphism $A/A^r_{k+1}\subset B/B^r_{k+1}$, and isomorphisms
$H_*(A/A^r_{k+1};\mathbb{Q})\cong H_*(B/B^r_{k+1};\mathbb{Q})$; and an isomorphism $(A^r_k/A^r_{k+1})\otimes \mathbb{Q}\cong (B^r_k/B^r_{k+1})\otimes
\mathbb{Q}$.
\item [2)] $\phi$ induces an epimorphism $H_2(A;\mathbb{Q})/<\Phi_n^r(A)>\to H_2(B;\mathbb{Q})/<\Phi_n^r(B)>$. 
\end{itemize}
\end{thm}

\begin{proof} First we prove that $1)\Rightarrow 2)$. This is trivial if $n=1$ so we assume $n>1$. Consider Stallings' exact sequence \cite[Section 7]{St}.
\begin{equation*}
\begin{CD}
 @>>>  H_2(A;\mathbb{Q}) @>\pi_A>> H_2(A/A^r_n;\mathbb{Q})  @>\partial_A>>
(A^r_n/A^r_{n+1})\otimes \mathbb{Q}   @>>> 0\\
&&    @VV \phi_* V        @VV(\phi_{n})_*V       @VV\phi_{n}V\\
 @>>>   H_2(B;\mathbb{Q})  @>\pi_B>>    H_2(B/B^r_n;\mathbb{Q})  @>\partial_B>>
(B^r_n/B^r_{n+1})\otimes
\mathbb{Q}    @>>> 0\\
\end{CD}
\end{equation*}
By $1)$, both $(\phi_n)_*$ and $\phi_n$ are isomorphisms. Noting that ker($\partial_A)= ~$image($\pi_A)\cong H_2(A;\mathbb{Q})/$ker($\pi_A)$, it follows that $\phi$ induces an isomorphism 
$$
\phi_*: H_2(A;\mathbb{Q})/ker(\pi_A)\cong H_2(B;\mathbb{Q})/ker(\pi_B).
$$
But ker($\pi_A)= ~<\Phi_n^r(A)>$, the subspace spanned by $\Phi_n^r(A)$ as may be seen by examining the diagram below.
\begin{equation*}
\begin{CD}
 0 @>>> \Phi_n^r(A) @>>> H_2(A;\mathbb{Z}) @>\pi_A^{\mathbb{Z}}>> H_2(A/A^r_n;\mathbb{Z}) \\
&&  @VVi_* V   @VVi_* V        @VV(i_n)_*V     \\
 0 @>>> ker(\pi_A) @>>> H_2(A;\mathbb{Q}) @>\pi_A>> H_2(A/A^r_n;\mathbb{Q}) \\ \\
\end{CD}
\end{equation*}
This completes the proof that $1)\Rightarrow 2)$.

Now we show $2)\Rightarrow 1)$. The proof follows the outline of Stallings' proof of his Rational Theorem (\cite[Theorem 7.3]{St}). We claim that the hypothesis $H_1(A;\mathbb{Q})\cong H_1(B;\mathbb{Q})$ establishes $1)$ for $k=1$. For note that $H_1(A;\mathbb{Q})\cong (A/A^r_2)\otimes \mathbb{Q}$. Furthermore, since $A/A^r_2$ is a torsion-free abelian group, 
$$
H_*(A/A^r_2;\mathbb{Q})\cong H_*((A/A^r_2)\otimes \mathbb{Q};\mathbb{Q})
$$
by \cite[Lemma 7.1]{St}. Proceeding inductively, we assume that for some fixed $k$, $1\leq k <n$, $\phi$ induces a monomorphism $A/A^r_{k}\subset B/B^r_{k}$, an isomorphism $H_*(A/A^r_{k};\mathbb{Q})\cong H_*(B/B^r_{k};\mathbb{Q})$, and an isomorphism
$(A^r_{k-1}/A^r_{k})\otimes \mathbb{Q}\cong (B^r_{k-1}/B^r_{k})\otimes \mathbb{Q}$ (similarly for all lesser values of $k$). For the inductive step we first show that $\phi$ induces an isomorphism
$(A^r_k/A^r_{k+1})\otimes \mathbb{Q}\cong (B^r_k/B^r_{k+1})\otimes
\mathbb{Q}$. For this we consider Stallings' exact sequence \cite[Section 7]{St}.
\begin{equation*}
\begin{CD}
 @>>>  H_2(A;\mathbb{Q}) @>\pi^k_A>> H_2(A/A^r_k;\mathbb{Q})  @>\partial_A>>
(A^r_k/A^r_{k+1})\otimes \mathbb{Q}   @>>> 0\\
&&    @VV \phi_* V        @VV(\phi_{k})_*V       @VV\phi_{k}V\\
 @>>>   H_2(B;\mathbb{Q})  @>\pi^k_B>>    H_2(B/B^r_k;\mathbb{Q})  @>\partial_B>>
(B^r_k/B^r_{k+1})\otimes
\mathbb{Q}    @>>> 0\\
\end{CD}
\end{equation*}
By the induction hypothesis, the middle map $(\phi_k)_*$ is an isomorphism. It follows immediately that $\phi_k$ is surjective. To establish injectivity of $\phi_k$, a diagram chase reveals that it suffices to show that $\phi$ induces an epimorphism $H_2(A;\mathbb{Q})/$ker($\pi_A^k)\to H_2(B;\mathbb{Q})/$ker($\pi_B^k$). By $2)$, $\phi$ induces an epimorphism 
$$
H_2(A;\mathbb{Q})/<\Phi_n^r(A)>\to H_2(B;\mathbb{Q})/<\Phi_n^r(B)>.
$$
We claim that $<\Phi_n^r(A)>= ~$ker($\pi_A^n$). For, if $x\in \Phi_n^r(A)$ then , by definition, $x\in ker(\pi^n_{\mathbb{Z}})$ below. Hence $i_*(x)\in ker(\pi^n_{\mathbb{Q}})$. Thus $\Phi_n^r(A)\subset ker(\pi_A^n)$. If $y\in ker(\pi_{\mathbb{Q}}^n)$ below then for some positive integer $m$, $my=i_*(x)$ for some $x\in ker(\pi^n_{\mathbb{Z}})=\Phi_n^r(A)$. Hence $y\in <\Phi_n^r(A)>$. Thus $<\Phi_n^r(A)>= ~$ker($\pi_A^n$).
\begin{equation*}
\begin{CD}
H_2(A;\mathbb{Z}) @>\pi^n_{\mathbb{Z}}>> H_2(A/A^r_n;\mathbb{Z})\\
    @VV i_* V        @VV i_* V       \\
H_2(A;\mathbb{Q})  @>\pi^n_{\mathbb{Q}}>>    H_2(A/A^r_n;\mathbb{Q}) \\
\end{CD}
\end{equation*}

The desired result now follows from the commutative diagram below since both horizontal maps $\pi$ are surjective.
\begin{equation*}
\begin{CD}
 H_2(A;\mathbb{Q})/ker(\pi_A^n) @>\pi>> H_2(A;\mathbb{Q})/ker(\pi_A^k)\\
    @VV \phi_* V        @VV\phi_*V       \\
H_2(B;\mathbb{Q})/ker(\pi_B^n) @>\pi>> H_2(B;\mathbb{Q})/ker(\pi_B^k)\\ 
\end{CD}
\end{equation*}
This completes our verification of the first part of the inductive step, namely that $(A^r_k/A^r_{k+1})\otimes \mathbb{Q}\cong (B^r_k/B^r_{k+1})\otimes
\mathbb{Q}$.

But this fact, together with the fact that the groups $A^r_k/A^r_{k+1}$ and $B^r_k/B^r_{k+1}$ are torsion-free, implies that $\phi$ induces an embedding $A^r_k/A^r_{k+1}\subset B^r_k/B^r_{k+1}$, which in turn implies that $\phi$ induces an embedding $A/A^r_{k+1}\subset B/B^r_{k+1}$.

It only remains to show that $H_*(A/A^r_{k+1};\mathbb{Q})\cong H_*(B/B^r_{k+1};\mathbb{Q})$. By our inductive assumption $H_*(A/A^r_{k};\mathbb{Q})\cong H_*(B/B^r_{k};\mathbb{Q})$. Moreover since $(A^r_k/A^r_{k+1})\otimes \mathbb{Q}\cong (B^r_k/B^r_{k+1})\otimes
\mathbb{Q}$, and since $A^r_k/A^r_{k+1}$ and $B^r_k/B^r_{k+1}$ are torsion-free abelian groups,
$$
H_*(A^r_k/A^r_{k+1};\mathbb{Q})\cong H_*(A^r_k/A^r_{k+1}\otimes \mathbb{Q};\mathbb{Q})
$$
$$
H_*(B^r_k/B^r_{k+1};\mathbb{Q})\cong H_*(B^r_k/B^r_{k+1}\otimes \mathbb{Q};\mathbb{Q})
$$
by \cite[Lemma 7.1]{St}. Thus 
$$
H_*(A^r_k/A^r_{k+1};\mathbb{Q})\cong H_*(B^r_k/B^r_{k+1};\mathbb{Q}).
$$
But the sequence
$$
0 \lra A^r_k/A^r_{k+1} \lra A/A^r_{k+1} \lra A/A^r_{k} \lra 0
$$
is a central extension so the result follows from \cite[Lemma 7.2]{St}(a Serre spectral sequence argument). This completes the proof of $2)\Rightarrow 1)$.

\end{proof}

\section{Applications}\label{applications}

\subsection{Algebraic Applications}

Suppose that $A$ is the subgroup generated by a set $\{a_1,\dots,a_m\}$ of elements of the group $B$.  When is $A$ a free group of rank $m$? When is $A$ free solvable? When is the image of $A$ in $B/B^{(n)}$ isomorphic to the free solvable group $F/F^{(n)}$?

\begin{prop}\label{freesolvable1} Suppose that $B$ is a finitely-related group and $A$ is the subgroup generated by $\mathcal{A}=\{a_i | i\in \mathcal{I}\}\subset B$. Suppose $\mathcal{A}$ is linearly independent in $H_1(B;\mathbb{Q})$ and suppose that $H_2(B;\mathbb{Q})= <\Phi^{(n-1)}(B)>$. Then $A/A^{(n)}$ is the free solvable group of derived length $n$ on $\mathcal{A}$, that is, if $F$ is the free group on $\mathcal{A}$ then the map $F\to A$ induces an isomorphism $F/F^{(n)}\cong A/A^{(n)}$. In particular $A$ maps onto the free solvable group on $\mathcal{A}$ of derived length $n$ and hence is not nilpotent if $m>1$. Moreover $A/A^{(n)}$ embeds in $B/B^{(n)}$.
\end{prop}

\begin{proof} By hypothesis $H_2(B;\mathbb{Q})/<\Phi_H^{(n-1)}(B)>=0$ and the map $\phi:F(\mathcal{A})\to B$ induces a monomorphism on $H_1(-;\mathbb{Q})$. Thus by Theorem~\ref{maintheorem}, $\phi$ induces a monomorphism $F/F^{(n)}\cong B/B^{(n)}$. This factors through the natural epimorphism $F/F^{(n)}\to A/A^{(n)}$ which is consequently an isomorphism. The other statements follow immediately.
\end{proof}

Examples of topological situations where the hypotheses of Proposition~\ref{freesolvable1} are satisfied are included in the next subsection.

\subsection{Topological Applications}

As previously discussed, Stallings' theorem has been instrumental in the study of link concordance. Recently, several other weaker equivalence relations on knots and links have been considered and found to be useful in understanding knot and link concordance ~\cite{COT}~\cite{COT2}~\cite{ConT1}~\cite{ConT2}~\cite{CT}~\cite{KT}~\cite{Ha2}~\cite{T}. These equivalence relations involved replacing the annuli in the definition of concordance by surfaces equipped with some extra structure. Below we show that our results generalize Stallings' results on link concordance to these more general equivalence relations. Moreover, recently, Harvey defined a rich new family of real-valued concordance invariants, $\rho_k(L)$, for a link $L$~\cite{Ha2}. She showed that these were actually invariants of some of these weaker equivalence relations and deduced new information about the Cochran-Orr-Teichner filtration of the classical disk-link concordance group. We are able to extend and refine her results. Details follow.

Recall that Stallings showed that concordant links have exteriors whose fundamental groups are isomorphic modulo any term of the lower central series. We can prove an analogue for the derived series and moreover use our Dwyer-type theorem to generalize his result to the following equivalence relation that is weaker than concordance.

\begin{defn}\label{ncobordant}Suppose that $L_0$ and $L_1$ are oriented, ordered, $m$-component links of circles in $S^3$. We say they are \emph{$(n)$-cobordant} if there exist compact oriented surfaces $\Sigma_i, 1\leq i\leq m$, properly and disjointly embedded in $S^3\times [0,1]$, restricting to yield $L_j$ on $S^3\times \{j\}$, $j=0,1$, and such that, for each $i$, for some set of circles $\{a_j,b_j\}$ representing a symplectic basis of curves for $\Sigma_i$, the image of each of the loops $\{a_j,b_j\}$ in $\pi_1((S^3\times [0,1])- \coprod \Sigma_i)\equiv \pi_1(E)$ is contained in $\pi_1(E)^{(n)}$ (use the unique `unlinked' normal vector field on $\Sigma_i$ to push off). $L_0$ is \emph{null $(n)$-bordant} if there are disjoint surfaces in $B^4$ as above whose boundaries form $L_0$.
\end{defn}

\begin{prop}\label{injectivity} Suppose that the $m$-component links $L_0$ and $L_1$ are $(n)$-cobordant via surfaces $\{\coprod \Sigma_i\}$ as above. Let $A$, $\overline{A}$ and $B$ denote the fundamental groups of their respective exteriors. Then both inclusion-induced maps $A\to B$ and $\overline{A}\to B$ satisfy the hypotheses of Theorem~\ref{main} for $n-1$. Consequently the rank of $A^{(n-1)}_H/A^{(n)}_H$ is the same as the rank of $\overline{A}^{(n-1)}_H/\overline{A}^{(n)}_H$. In addition the inclusion maps induce monomorphisms $A/A^{(n)}_H\hookrightarrow B/B^{(n)}_H$ and $\overline{A}/\overline{A}^{(n)}_H\hookrightarrow B/B^{(n)}_H$. If $L_0$ is null $(n)$-cobordant then rank of $A^{(n-1)}_H/A^{(n)}_H$ is the same as that of the $m$-component trivial link, namely $m-1$ (for $n\geq 2)$, and the set of meridians viewed as a subset of either $A$ or $B$ satisfies the hypotheses of Proposition~\ref{freesolvable1}.
\end{prop}

\begin{proof} We use the notation of Definition~\ref{ncobordant} and Proposition~\ref{injectivity}. By hypothesis, for each $1\leq i\leq m$ there exist symplectic bases of circles $\{a_{ij},b_{ij}\}$ for $\Sigma_i$ whose push-offs $\{a^+_{ij},b^+_{ij}\}$ into $E$ lie in $\pi_1(E)^{(n)}= B^{(n)}$. The key observation is that $H_2(E;\mathbb{Z})$ is generated by the tori $\{a^+_{ij}\times S^1_i, b^+_{ij}\times S^1_i\}$ where $S^1_i$ is a fiber of the normal circle bundle to $\Sigma_i$, together with the $m$ tori $L_0\times S^1_i$ that live in $S^3-L_0$. Thus the cokernel of the map $H_2(A;\mathbb{Z})\to H_2(B;\mathbb{Z})$ is generated by the former collections. Since $[a^+_{ij}]\in B^{(n)}$, $a^+_{ij}$ bounds a $B^{(n-1)}$-surface $S_{ij}$ mapped into $E$ (that is $\pi_1(S_{ij})\subset B^{(n-1)}$). If we cut open the torus $a^+_{ij}\times S^1_i$ along $a^+_{ij}$ and adjoin two oppositely oriented copies of $S_{ij}$, we obtain a (mapped in) surface that is homologous to $a^+_{ij}\times S^1_i$ and is also a $B^{(n-1)}$-surface. Similarly for the tori $b^+_{ij}\times S^1_i$. Therefore $A\to B$ satisfies the hypotheses of Theorem~\ref{main} for $n-1$. Thus this map induces a monomorphism $A/A^{(n)}_H\hookrightarrow B/B^{(n)}_H$ and the ranks of $A^{(n-1)}_H/A^{(n)}_H$ and $B^{(n-1)}_H/B^{(n)}_H$ over their respective rings are equal. Symmetrically, the same is true for $\overline{A}\to B$. The first part of the theorem follows immediately.

If $L_0$ is null-$(n)$-cobordant then $\overline{A}=F$ the free group of rank $m$. It is known that the rank of $F^{(n-1)}/F^{(n)}$ is $m$ if $n=1$ and $m-1$ if $n\geq 2$ ~\cite[Lemma 2.12]{COT}. Moreover, since the longitudes of the components of $L_0$ co-bound the $\Sigma_i$ with the longitudes of the trivial link, which are trivial, the longitudes of $L_0$ map into $B^{(n+1)}$. By the first part of the theorem, this implies they lie in $A^{(n)}$ and hence bound (immersed) $A^{(n-1)}$-surfaces $S_i$ in $S^3-L_0$. On the other hand, $H_2(S^3-L_0;\mathbb{Z})$ is generated by the $m$ tori $L_0\times S^1_i$ (the boundaries of the regular neighborhoods of the components of $L_0$). The $S_i$ can be used to surger these tori and thus showing that $H_2(A;\mathbb{Z})$ is generated by $A^{(n-1)}$-surfaces. Hence the meridional set in $A$ satisfies the hypotheses of Proposition~\ref{freesolvable1}. Above we saw that the cokernel of the map $H_2(A;\mathbb{Z})\to H_2(B;\mathbb{Z})$ was generated by $B^{(n-1)}$-surfaces. Combining these two facts, we see that $H_2(B;\mathbb{Z})$ is generated by $B^{(n-1)}$-surfaces. Hence the set of meridians of $L_0$ viewed as a subset of $B$ satisfies the hypotheses of Proposition~\ref{freesolvable1}.
\end{proof}

In ~\cite{Ha2}, Harvey defined a new family of real-valued invariants, $\rho_k$, $k\geq 0$, of closed odd-dimensional manifolds using the torsion-free derived series and a higher-order signature defect, the Cheeger-Gromov von Neumann $\rho$-invariant. Let $M_L$ be the closed 3-manifold associated to $(S^3,L)$ by performing $0$-framed Dehn surgery on $S^3$ along the components of $L$ ~\cite{R}. Let $G= \pi_1(M_L)$. Harvey associates to $L$ the pair $(M_L, \phi_k:G\to G/G^{(k+1)}_H)$ where $G^{(k+1)}_H$ is the ($k+1$)-st term of Harvey's torsion-free derived series. Then she defines $\rho_k(L)\equiv \rho(M_L,\phi_k)$ where the latter is the Cheeger-Gromov von Neumann $\rho$-invariant ~\cite{ChGr1}. Actually Harvey defines and establishes these invariants independent of the work of Cheeger and Gromov, but observes that they coincide with the invariants of Cheeger-Gromov. Harvey established that these link invariants were concordance invariants by showing the manifold invariants were rational homology cobordism invariants. In particular, all of these invariants vanish for links concordant to the trivial link, called \emph{slice links}, that is links whose components bound disjoint embedded disks in the $4$-ball. But she went on to show that the $\rho_k$ actually respected some even weaker equivalence relations ~\cite[Theorem 6.4]{Ha2}. In particular she showed that $\rho_n(L)$ vanishes for links in $\mathcal{F}_{(n+1)}$, the set of $(n+1)$-solvable links, that is the ($n+1$)-st term of the filtration of the link concordance group defined in ~\cite[Section 8]{COT} (these are reviewed herein). This class is much larger than that of slice links.  We improve on her result. Recall that an $m$-component link $L$ is a $\textbf{finite E-link}$ if $\pi_1(S^3-L)$ admits a homomorphism to a finite $E$-group of rank $m$ under which the longitudes map trivially (a \textbf{finite E-group} is one that is the fundamental group of a finite 2-complex with $H_1\cong \mathbb{Z}^m$ and $H_2\cong 0$). Boundary links, homology boundary links, fusions of boundary links and sublinks of homology boundary links are all finite $E$-links \cite[p.641-644]{C2}. Conjecturally, the class of finite E-links is the same as the class of links with vanishing Milnor's $\ov{\mu}$-invariants. Our generalization is the following.

\begin{thm}\label{n.5solvable} $\rho_n(L)$ vanishes for all finite $E$-links in $\mathcal{F}^{\mathbb{Q}}_{(n.5)}$, the set of all rationally (n.5)-solvable links (see ~\cite[Section 4, Section 8]{COT}).
\end{thm}

This theorem improves on Harvey's ~\cite[Theorem 6.4]{Ha2} in two ways. Firstly and primarily, it improves the $(n+1)$ in her result to what should be the optimal result $(n.5)$ (although only for $E$-links). Being an $E$-link ensures an extra rank condition that is hidden in Harvey's proof since it is implied by $(n+1)$-solvability. Our Theorem then allows for a corresponding sharpening of Harvey's main application of the above ~\cite[Theorem 6.8]{Ha2}.

\begin{thm}\label{infgenerated} In the category of m-component ordered oriented string links ($m>1$), each of the quotients $\mathcal{F}_{(n)}/\mathcal{F}^{\mathbb{Q}}_{(n.5)}$ contains a subgroup, consisting entirely of boundary links, whose abelianization has infinite $\mathbb{Q}$-rank.
\end{thm}

We remark that Harvey's additivity result for her $\rho_k$ for boundary string links should hold for any additively closed subset of $E$-links, such as the set of homology boundary links with a fixed ``\emph{pattern}'' and so the word ``boundary links'' in the above theorem should be able to be replaced by any such set.

Secondly, our Theorem~\ref{n.5solvable} improves on Harvey's version by proving the theorem for so-called ``rational'' solvability. The definition of the latter is reviewed below. This may seem like a technical advance. However it points the way to certain further improvements in other results in the literature that we will postpone to another paper. Namely, there is a well-studied geometric notion that approximates $n$-solvability of links that has to with generalizing the annuli in the definition of concordance to certain $2$-complexes called \emph{symmetric gropes}. We shall not review these terms here, but in future paper we will show that the torsion-free derived series suggests beautiful generalizations of gropes, that we call \emph{rational homology gropes} and these generalizations are the proper geometric approximation to the algebraic notion of rational n-solvability.

We briefly review the definitions of the Cheeger-Gromov von Neumann $\rho-$invariant (only in the cases that we need here) and $(n)-$solvability. More general definitions are to be found in ~\cite{ChGr1}~\cite[Section 3]{Ha2}~\cite[Section 5]{COT}.

\begin{defn}\label{def:rhoinvariant} Suppose $M$ is a closed, oriented $3-$manifold and $\phi:\pi_1(M)\to \G$is a homomorphism to a poly-(torsion-free-abelian) group. Suppose also that $\phi$ extends to $\psi:\pi_1(W)\to \G$ where $W$ is a compact, oriented $4$-manifold whose boundary is $M$. Then $\rho(M, \phi)$ is given by $\sigma^{(2)}_{\Gamma}(W)-\sigma(W)$, where $\sigma^{(2)}_{\Gamma}(W)$ is the von Neumann signature of the equivariant intersection form $\lambda_{\Gamma}$ on $H_2(W;\mathcal{K}\Gamma)$ and $\sigma(W)$ is the signature of the usual intersection form on $H_2(W;\mathbb{Q})$.
\end{defn}

\begin{defn}\label{def:ratnsolvable}(see ~\cite[Section 4]{COT}) A connected, closed, oriented $3$-manifold $M$ is $\textbf{rationally (n)-solvable}$ if there exists a compact, connected, oriented $4$-manifold $W$ such $\partial(W)=M$ and
\begin{itemize}
\item The inclusion map induces an isomorphism $j_*:H_1(M;\mathbb{Q})\to H_1(W;\mathbb{Q})$,
\item $H_2(W;\mathbb{Q})$ admits a basis $\{[L_i],[D_i]; 1\leq i\leq r\}$ consisting of connected, oriented, $\pi_1(W)^{(n)}$-surfaces $\{L_i, D_i\}$ (that is $\pi_1(L_i)$ and $\pi_1(D_i)$ are contained in $\pi_1(W)^{(n)}$), whose equivariant intersection numbers with coefficients in $\bbq[ \pi_1(W)/\pi_1(W)^{(n)}]$ are as follows $L_i\cd L_j =0$ and $D_i\cd D_j =0$ if $i\neq j$ and $L_i\cd D_i= 1$. 
\end{itemize}
In this case the $L_i$ are said to generate an $\textbf{(n)-Lagrangian}$ for $W$ and the $D_i$ are said to generate its $\textbf{(n)-duals}$. In this case we say that $M$ is $\textbf{rationally}$ $\textbf{(n)-solvable via W}$.
A $3$-manifold $M$ is $\textbf{rationally (n.5)-solvable}$ if it satisfies the above and in addition:
\begin{itemize}
\item The $L_i$ are $\pi_1(W)^{(n+1)}$-surfaces which are then said to constitute an 
$\textbf{(n+1)}$-$\textbf{Lagrangian}$ for $W$.
\end{itemize}
\end{defn}

\begin{defn}\label{def:linkratnsolvable} A link $L$ in $S^3$ is said to be $\textbf{rationally (n)-solvable}$ (respectively $\textbf{rationally (n.5)-solvable}$) if the zero-framed surgery $M_L$ is rationally $(n)-$solvable (respectively rationally $(n.5)-$solvable) as above. The set of (concordance classes) of such links is denoted $\mathcal{F}^{\mathbb{Q}}_{(n)}$ (respectively $\mathcal{F}^{\mathbb{Q}}_{(n.5)}$).
\end{defn}

In ~\cite[Sections 4 and 8]{COT} these notions and also the notions of $\textbf{(n)-solvable}$ and $\textbf{(n.5)-solvable}$ were defined. The latter are the same as the above except that $\mathbb{Q}$ is replaced by $\mathbb{Z}$, $W$ is required to be spin and the equivariant self-intersection form is also considered. Links satisfying these stronger conditions are said to lie in $\mathcal{F}_{(n)}$ and $\mathcal{F}_{(n.5)}$ respectively. Note that a link that is $(n)-$solvable is certainly rationally $(n)-$solvable, and that rationally $(n.5)-$solvable links are certainly rationally $(k)-$solvable for any integer or half-integer $k\leq n.5$. It is easy to see that any slice link is $(n)-$solvable for all $n$ and if two links are concordant then one is $(n)-$solvable if and only if the other is also.

Note that a $\pi_1(W)^{(n+1)}$-surface of $W$ (sometimes called an $(n+1)-$surface), lifts to the regular $\G$ covering space of $W$ corresponding to a map $\psi:\pi_1(W)\to \G$ as long as $\G^{(n+1)}=1$ since then $\pi_1(W)^{(n+1)}\subset \text{ker}(\psi)$. Thus an $(n.5)$-Lagrangian $L$ lifts to generate a $\mathbb{Z}\G-$submodule of $H_2(W;\mathbb{Z}\G)$. By defintion the $\mathbb{Z}\G-$equivariant intersection form (on $H_2(W;\mathbb{Z}\G)$) is identically zero on this submodule. The same holds true for $H_2(W;\mathcal{K}\G)$. Thus if the intersection form is nonsingular and if this submodule has one-half rank, the von Neumann signature of the equivariant intersection form will be zero.

Theorem~\ref{n.5solvable} and Theorem~\ref{infgenerated} are consequences of the more general theorem below (and of Harvey's previous work). This theorem is a basic and important result in its own right. It generalizes ~\cite[Theorem 4.1]{COT} where this precise theorem is stated under the assumption that $\beta_1(M)=1$. The proof in ~\cite{COT} does \emph{not} apply to the general case. It also generalizes \cite[Theorem 6.4]{Ha2}.

\begin{thm}\label{thm:rho=0} Let $\Gamma$ be a PTFA group such that $\G^{(n+1)}=0$. Let $M$ be a closed, connected, oriented $3$-manifold equipped with a non-trivial coefficient system $\phi:\pi_1(M)\to \Gamma$. Suppose $\rank_{\mathbb{Z}\Gamma}(H_1(M;\mathbb{Z}\G))= \beta_1(M)-1$. If $M$ is
rationally
$(n.5)$-solvable via a
$4$-manifold $W$ over which $\phi$ extends, then
$\rho(M, \phi) = 0$.

Moreover, if $W$ is a rational ($n+1$)-solution then the above rank condition is automatically satisfied.
\end{thm}

The following corollary generalizes ~\cite[Theorem 6.4]{Ha2} where the hypothesis is that $M$ is $(n+1)-$solvable. We weaken this hypothesis to $(n.5)-$solvability but require an extra rank requirement.

\begin{cor}\label{cor:n.5solvable}
 If $M$ is rationally $(n.5)$-solvable and $$\text{rank}_{\mathbb{Z}[\pi/\pi_{\sss H}^{\sss (n+1)}]}
  \frac{\pi_{\sss H}^{(n+1)}}{\pi_{\sss H}^{(n+2)}}=\beta_{1}(M)-1$$
 where $\pi=\pi_{1}(M)$ then $\rho_{n}(M)=0.$ 
\end{cor}

\begin{proof}[Proof that Corollary~\ref{cor:n.5solvable} implies Theorem~\ref{n.5solvable}] Suppose $L$ is a finite $E$-link that is rationally $(n.5)$-solvable. Then, by definition, $M_L$ is rationally $(n.5)$-solvable. Taking $M=M_L$ and applying Corollary~\ref{cor:n.5solvable}, to conclude that $\rho_n(L)$ vanishes we need only verify that for an $E$-link $L$, $M_L$ satisfies the rank hypothesis of Corollary~\ref{cor:n.5solvable}. Let $\pi=\pi_1(M_L)$. By the definition of an $\textbf{E}$-link, there is a map $\phi:\pi\to E$ that is an isomorphism on $H_1(-;\mathbb{Z})$ and an epimorphism on $H_2(-;\mathbb{Z})$. Then, by ~\cite[Theorem 4.1]{CH1}, for any $n$,
$$
\text{rank}_{\mathbb{Z}[\pi/\pi_{\sss H}^{\sss (n+1)}]}
  \frac{\pi_{\sss H}^{(n+1)}}{\pi_{\sss H}^{(n+2)}}= \text{rank}_{\mathbb{Z}[E/E_{\sss H}^{\sss (n+1)}]}
  \frac{E_{\sss H}^{(n+1)}}{E_{\sss H}^{(n+2)}}.
  $$
By Proposition~\ref{observation2}, the latter expression is the same as rank($H_1(E,\mathbb{Z}[E/E^{(n+1)}_H]$). But by ~\cite[Lemma 2.12]{COT} the latter rank is precisely $\beta_{1}(E)-1= \beta_{1}(M_L)-1$ except in the degenerate case that $E= E^{(n+1)}_H$. This case cannot occur here since the meridional map $F\to \pi\to E$ is also homologically $2$-connected and so $F/F^{(n+1)}$ embeds in $E/E^{(n+1)}_H$ by ~\cite[Theorem 4.1, Proposition 2.4]{CH1}. 
\end{proof}

\begin{proof}[Proof of Theorem~\ref{infgenerated}] For any $n$ Harvey produced in ~\cite[Theorem 6.8]{Ha2} an infinite set of $(n)$-solvable boundary string links whose set of $\rho_n$ ( real numbers) was $\mathbb{Q}$-linearly independent. She also proved that $\rho_n$ was additive on the subgroup of boundary string links ~\cite[Corollary 6.7]{Ha2}. If any linear combination of these links were $(n.5)$-solvable, it would contradict our Theorem~\ref{n.5solvable}.
\end{proof}

The Proposition below was proved by Harvey under the slightly stronger hypothesis of $(n)-$solvability ~\cite[Theorem 6.4]{Ha2} (as opposed to \emph{rational} $(n)-$solvability). The interesting thing about our proof of the result below is that does not use the hypothesis that $M$ and $W$ are manifolds and it does not use anything about the intersection form. Rather it just uses the fact that $H_2(W;\mathbb{Q})$ has a basis of $(n)$-surfaces (and the main theorem of the current paper).

\begin{prop}\label{prop:nsolvable} If $M$ is rationally $(n)-solvable$ via $W$ then, letting $\pi=\pi_{1}(M)$ and $B=\pi_{1}(W)$, the inclusion $j:M \rightarrow W$ induces a monomorphism
$$j_{n+1}: \frac{\pi}{\pi_{\sss H}^{(n+1)}} \hookrightarrow  \frac{B}{B_{\sss H}^{(n+1)}}.$$
and if $n>0$,
$$\text{rank}_{\mathbb{Z}[\pi/(\pi)_{\sss H}^{\sss (n)}]}
  \frac{(\pi)_{\sss H}^{(n)}}{(\pi)_{\sss H}^{(n+1)}}=\beta_{1}(M)-1.
  $$
  
More generally, if $\psi:\pi_1(W)=B\to \G$ is any non-trivial PTFA coefficient system with $\G^{(n)}=1$, then
$$
\text{rank}_{\mathbb{Z}\G}H_1(M;\mathbb{Z}\G)=\beta_1(M)-1.
$$  
\end{prop}

\begin{cor}\label{cor:n+1solvable}(compare ~\cite[Theorem 6.4]{Ha2})
 If $M$ is rationally $(n+1)$-solvable and $\pi=\pi_1(M)$ then
 \begin{itemize}
 \item $\text{rank}_{\mathbb{Z}[\pi/\pi_{\sss H}^{\sss (n+1)}]}
  \frac{\pi_{\sss H}^{(n+1)}}{\pi_{\sss H}^{(n+2)}}=\beta_{1}(M)-1$
 and
 \item $\rho_{n}(M)=0.$ 
 \end{itemize}
\end{cor}

\begin{proof}[Proof of Proposition~\ref{prop:nsolvable}] We apply our main theorem (Theorem~\ref{main}). Since $W$ is a rational $n$-solution for $M$, $H_{2}(W;\mathbb{Q})$ has a basis of $B^{\sss (n)}$-surfaces so it certainly has a basis of $B_{\sss H}^{\sss (n)}$-surfaces since $B^{\sss (n)} \subset B^{\sss (n)}_{\sss H}$. The first statement of the Proposition then follows immediately.

For the second part, consider the meridional map $F\to \pi\to B$ and observe that Theorem~\ref{main} applies to it also. It follows that

$$\text{rank}_{\mathbb{Z}[\pi/(\pi)_{\sss H}^{\sss (n)}]}
  \frac{(\pi)_{\sss H}^{(n)}}{(\pi)_{\sss H}^{(n+1)}}\cong \text{rank}_{\mathbb{Z}[B/B_{\sss H}^{\sss (n)}]}
  \frac{B_{\sss H}^{(n)}}{B_{\sss H}^{(n+1)}}
  $$
  and 
  $$
  \text{rank}_{\mathbb{Z}[F/F_{\sss H}^{\sss (n)}]}
  \frac{F_{\sss H}^{(n)}}{F_{\sss H}^{(n+1)}}\cong \text{rank}_{\mathbb{Z}[B/B_{\sss H}^{\sss (n)}]}
  \frac{B_{\sss H}^{(n)}}{B_{\sss H}^{(n+1)}}.$$
  
  But $F^{(n)}_H=F^{(n)}$ by \cite[Proposition 2.3]{Ha2} and the rank of $F^{(n)}/F^{(n+1)}$ ($n>0$) is $\b_1(M)-1$
  by an easy Euler characteristic argument \cite[Lemma 2.12]{COT}.
  
For the last claim, consider $F\to \pi\to B\overset{\psi}\to \G$ and observe that Proposition~\ref{2-connected} applies to $F\to B$ and $\pi\to B$, since $H_2(W;\mathbb{Q})$ has a basis of ker($\psi$)-surfaces since $B^{(n)}\subset \text{ker}(\psi)$. Thus
$$
H_1(F;\mathcal{K}\G)\cong H_1(M;\mathcal{K}\G)\cong H_1(B;\mathcal{K}\G).
$$
As above, the rank of the first of these three is known to be $\beta_1(M)-1$ \cite[Lemma 2.12]{COT}.

\end{proof}

\begin{proof}[Proof that Theorem~\ref{thm:rho=0} and Proposition~\ref{prop:nsolvable} imply Corollary~\ref{cor:n.5solvable}] If $M$ is rationally $(n.5)$-solvable via $W$ then it is rationally $(n)$-solvable via $W$. Thus Proposition~\ref{prop:nsolvable} applies. Let $\pi=\pi_{1}(M)$ and $B=\pi_{1}(W)$. Since $j_{n+1}$ is injective, 
$$
\rho_{n}(M,\pi:\pi \rightarrow \pi/\pi_{\sss H}^{\sss (n+1)}))= 
\rho(M,j_{n+1}\circ \pi)
$$
by the $\G$-induction property of von Neumman $\rho$-invariants. Letting  $\psi$ be the canonical map $B \rightarrow B/B_{\sss H}^{\sss (n+1)}$ and letting $\phi= j_{n+1}\circ \pi= \psi\circ j$, we see that $\rho_n(M)= \rho(M,\phi)$. Since $\phi$ extends over $W$ by $\psi$, we may apply Theorem~\ref{thm:rho=0}
with $\Gamma=B/B_{\sss H}^{\sss (n+1)}$ and note that $\Gamma^{(n+1)}=\{e\}$ since $B^{\sss (n+1)}\subset \pi_1(W)_{\sss H}^{\sss (n+1)}$). We need only verify that the rank hypothesis of Corollary~\ref{cor:n.5solvable} implies the rank hypothesis of Theorem~\ref{thm:rho=0}. Let $\Gamma '= \pi_1(M)/\pi_1(M)^{(n+1)}_H$. Then, by Proposition~\ref{observation2}, the rank hypothesis of Corollary~\ref{cor:n.5solvable} is equivalent to the fact that $\mathcal{K}\Gamma '$-rank of $H_1(\pi_1(M);\mathcal{K}\Gamma ')$ is $\beta_1(M)-1$. Since $j_{n+1}: \Gamma '\to\G$ is a monomorphism, by \cite[Lemma 4.2]{CH1}, this rank is the same as the $\mathcal{K}\Gamma$-rank of $H_1(\pi_1(M);\mathcal{K}\Gamma)$ associated to the coefficient system $\phi$. But this is precisely the rank hypothesis of Theorem~\ref{thm:rho=0} in the case that $\Gamma=B/B_{\sss H}^{\sss (n+1)}$.
\end{proof}

\begin{proof}[Proof of Theorem~\ref{thm:rho=0}] Note that the final claim of the theorem was already established by the last part of Proposition~\ref{prop:nsolvable} applied to the $(n+1)$-solution $W$. 

Suppose $M$ is rationally $(n.5)$-solvable via $W$ such that the coefficient system extends to $\psi:\pi_1(W)\to \Gamma$. Then $\rho(M, \phi)$ is given by $\sigma^{(2)}_{\Gamma}(W)-\sigma(W)$, where $\sigma^{(2)}_{\Gamma}(W)$ by Definition~\ref{def:rhoinvariant}. Let $\tilde{I}$ denote the image of the map
$$
H_2(\partial W;\mathcal{K}\G)\overset{j_*}{\lra}H_2(W;\mathcal{K}\G).
$$
Since the sequence
$$
0\lra \tilde {I} \lra H_2(W;\mathcal{K}\G)\lra H_2(W;\mathcal{K}\G)/\tilde {I}\lra 0
$$
is split exact (since all $\mathcal{K}\Gamma$-modules are free) it follows that
$$
H_2(W;\mathcal{K}\G)\cong \tilde{I}\oplus (H_2(W;\mathcal{K}\G)/\tilde {I}).
$$
Since the intersection form
$$
\lambda: H_2(W;\mathcal{K}\G)\to Hom(H_2(W;\mathcal{K}\G), \mathcal{K}\G)
$$
is the composition of $\pi_*$ (below)
$$
H_2(\partial W;\mathcal{K}\G)\overset{j_*}{\lra}H_2(W;\mathcal{K}\G)\overset{\pi_*}{\lra} H_2(W,\partial W;\mathcal{K}\G)
$$ 
followed by the Poincar\'{e} Duality and the Kronecker map (both of the latter are isomorphisms) the kernel of $\lambda$ is the kernel of $\pi_*$ which is precisely $\tilde{I}$. Hence $\lambda$ induces a \emph{nonsingular} intersection form
$$
\tilde \lambda: H_2(W;\mathcal{K}\G)/\tilde{I}\to Hom((H_2(W;\mathcal{K}\G)/\tilde {I}),\mathcal{K}\G)
$$ 
and, with respect to the direct sum decomposition above, $\lambda$ is the direct sum of $\tilde\lambda$ and the zero form on $\tilde{I}$. Hence
$$
\sigma^{(2)}_{\Gamma}(W)= \sigma^{(2)}_{\Gamma}(\lambda) = \sigma^{(2)}_{\Gamma}(\tilde\lambda).
$$ 
Therefore, since $\tilde\lambda$ is nonsingular, we need to show that there is a one-half rank summand of $H_2(W;\mathcal{K}\G)/\tilde {I}$ on which $\tilde\lambda$ vanishes.

Towards this end, first recall that
$$
H_1(M;\mathbb{Q})\overset{j_*}{\lra} H_1(W;\mathbb{Q}),
$$
is an isomorphism 
$$
H_3(W;\mathbb{Q})\cong H^1(W,M;\mathbb{Q})=0, 
$$
and
$$
H_3(W,M;\mathbb{Q})\cong H^1(W;\mathbb{Q})\cong H^1(M;\mathbb{Q}) \cong H_2(M;\mathbb{Q})
$$
so 
$$
H_3(W,M;\mathbb{Q})\overset{\partial_*}{\lra}H_2(M;\mathbb{Q})
$$
is a monomorphism between vector spaces of the same rank, hence an isomorphism. It follows that

$$
H_2(W;\mathbb{Q})\overset{\pi_*}{\lra} H_2(W,M;\mathbb{Q})
$$
is an isomorphism. Therefore $\rank_\mathbb{Q}(H_2(W;\mathbb{Q}))= \rank_\mathbb{Q}(H_2(W,M;\mathbb{Q}))= \beta_2(W)$. Let $\Sigma_1,\dots, \Sigma_r$ be $\pi_1(W)^{(n+1)}$-surfaces of $W$ representing a rational $(n+1)$-lagrangian for $W$. By hypothesis, $\Sigma_1,\dots, \Sigma_r$ is the basis of a one-half-rank $\mathbb{Q}$-vector space $L$ in
$H_2(W;\mathbb{Q})$. Then $\rank
_\mathbb{Q}(L)= \rank
_\mathbb{Q}(\pi_*(L))= (1/2)\beta_2(W)$. Hence
$$
\rank _\mathbb{Q}(H_2(W,M;\mathbb{Q})/\pi_*(L))= (1/2)\beta_2(W).
$$
Since $\G^{(n+1)}=0$,
$\psi$ factors through the quotient
$\pi_1(W)/\pi_1(W)^{(n+1)}$ so $\pi_1(\Sigma_i)\subset ker\psi_i$. Using this let $\tilde{L}$ be the submodule
generated by $\tilde{\Sigma_1},\dots, \tilde{\Sigma_r}$ in
$H_2(W;\mathcal{K}\Gamma)/\tilde{I}$. By naturality, the intersection form with $\mathcal{K}\Gamma$ coefficients,
$\lambda$ vanishes on $\tilde{L}$. Thus $\tilde{L}$ is a free summand of $H_2(W;\mathcal{K}\Gamma)/\tilde{I}$ that is isomorphic to its image $\pi_*(\tilde{L})$ in $H_2(W,M;\mathcal{K}\G)$. Therefore, to conclude that $\sigma^{(2)}_{\Gamma}(\tilde\lambda)=0$, it suffices to show that 
$$
\rank_{\mathcal{K}\Gamma}( 
\tilde{L})\geq (1/2)(\rank_{\mathcal{K}\Gamma}(H_2(W;\mathcal{K}\Gamma)/\tilde{I}). 
$$
Let $b_i(M)= \rank_{\mathcal{K}\G}(H_i(M;\mathcal{K} \G))$, and $b_i(W)= \rank_{\mathcal{K}\G}(H_i(W;\mathcal{K} \G))$. Note  also $b_i(W)=\rank(H_{4-i}(W,M;\mathcal{K} \G))$. Since $W$ is a topological 4-manifold with non-empty boundary, it has the homotopy type of $3$-dimensional CW complex. Since $H_3(W;\mathbb{Q})\cong H^1(W,\partial W;\mathbb{Q})=0$, the boundary homomorphism
$\partial_3:C_3(W)\to C_2(W)$ is injective. Let
$C_{\ast}(W;\mathbb{Q}\G)$ be the corresponding $\mathbb{Q}\G$ chain complex free on the cells of
$W$. By Strebel's ~\cite[p.305]{Str},
$\tilde{\partial}_3 :C_3(W;\mathbb{Q}\G)\to C_2(W;\mathbb{Q}\G)$ is injective so
$H_3(W;\mathbb{Q}\G)=0$. It follows that both $H_3(W;\mathcal{K}\G)=0$ and $H_1(W,M;\mathcal{K}\G)=0$ so
$b_3=0$ . Thus
$$
H_1(M;\mathcal{K}\G)\overset{j_*}{\lra} H_1(W;\mathcal{K}\G),
$$
is an epimorphism whose kernel, $K$, has rank $b_1(M)- b_1(W)$, and since
$$
0\lra H_2(W;\mathcal{K}\G)/\tilde{I}\overset{\pi_*}{\lra} H_2(W,M;\mathcal{K}\G)\overset{\p_*}{\lra}K\lra 0
$$
is exact, 
$$
\rank_{\mathcal{K}\Gamma}(H_2(W;\mathcal{K}\Gamma)/\tilde{I})=b_2(W)+b_1(W)-b_1(M).
$$
Therefore our goal translates to showing that
$$
\rank_{\mathcal{K}\Gamma}( 
\tilde{L})\geq 1/2(b_2(W)+b_1(W)-b_1(M)).
$$
Towards this goal, apply Corollary~\ref{cor:independence} setting $\mathcal{C}_*= C_*(W,M;\mathbb{Q}\Gamma)$, $\ov{\mathcal{C}}_*=C_*(W,M;\mathbb{Q})$, $p=2$ and $\{x_s\}= \{\Sigma_1,\dots, \Sigma_r\}$ to conclude that
$$
\rank_{\mathcal{K}\Gamma}(H_2(W,M;\mathcal{K}\Gamma)/\pi_*(\tilde{L}))\leq \rank _\mathbb{Q}(H_2(W,M;\mathbb{Q})/\pi_*(L)),
$$
and thus that
$$
\rank_{\mathcal{K}\Gamma}(H_2(W,M;\mathcal{K}\Gamma)/\pi_*(\tilde{L}))\leq (1/2)\beta_2(W).
$$
We conclude that
$$
\rank_{\mathcal{K}\Gamma}( 
\tilde{L})\geq b_2(W) - (1/2)\beta_2(W).
$$
Hence, we will have achieved our goal above if we can show that
$$
b_2 - (1/2)\beta_2(W)\geq 1/2(b_2(W)+b_1(W)-b_1(M)).
$$
By our hypothesis, $b_1(M)= \beta_1(M)-1= \beta_1(W)-1$ so the last inequality is equivalent to
$$
b_2(W)-b_1(W)\geq \beta_2(W)-\beta_1(W)+1= \chi (W).
$$
But this is surely true since in fact the Euler characteristic of $W$ can be computed using $\mathcal{K}\G$-coefficients so $\chi (W)= b_2(W)-b_1(W)$ since $b_4(W)=b_3(W)=b_0(W)=0$. Actually for $b_0(W)=0$ we need that $\beta_1(W)=\beta_1(M)\geq 1$ to ensure that the coefficient system $\psi$ is non-trivial (~\cite[Prop.2.9]{COT}). In the case that $\beta_1(W)=\beta_1(M)=0$, both $\phi$ and $\psi$ are trivial since $\G$ is poly-(torsion-free-abelian). This case was excluded by the hypotheses. However note that in this degenerate case, the von Neumann signature and the ordinary signature are identical so $\rho(M,\phi)=0$ automatically.
\end{proof}

\section{Homological Localization}\label{HomologicalLocalization}
Recall that in ~\cite{CH1} the authors constructed a rational homological localization, $G\to \wt G$,  called \emph{the torsion-free-solvable completion}. This means that a rational homology equivalence $A\to B$ induces an isomorphism $\wt A\to \wt B$. In the context of rational homological localization, it was suggested that this could be viewed as an analogue of the Malcev completion, $G\otimes \mathbb{Q}$, of a group $G$ wherein one replaces the lower central series by the torsion-free derived series. Recall that $G\otimes \mathbb{Q}$ is defined to be the inverse limit of certain $n$-torsion-free-nilpotent groups $G/G_n\otimes \mathbb{Q}$. Similarly $\wt G$ was defined to be the inverse limit of a certain tower of \emph{$n$-torsion-free-solvable groups} $\wt G_n$. Here we can use our version of Dwyer's theorem for the torsion-free derived series to prove that $\wt G_n$ is functorially preserved by a larger class of maps (than rational homology equivalences).

To state the result, we recall a few definitions. The reader is referred to ~\cite{CH1} for more detail.

\begin{defn}\label{deftorsionfreesolv} A group $A$ is \textbf{$n$-torsion-free-solvable} if $A^{\sss (n)}_{\sss H}=0$ .
\end{defn}

\begin{defn}\label{torsionfreetower}(compare \cite[Section 12]{B}) A collection of groups $A_n$, $ n\geq 0$, and group homomorphisms $f_n, \pi_n$ , $n\geq 0$, as below:
$$
\begin{array}{cccccc}
A     &\overset{f_n}{\lra}    &  A_n   &\overset{\pi_n}{\lra}
& A_{n-1}\\
\end{array}
$$
compatible in the sense that $f_{n-1}=\pi_n\circ f_n$,  is a \textbf{torsion-free-solvable tower for A} if, for each $n$, $A_n$ is $n$-torsion-free-solvable and the kernel of $\pi_n$ is contained in $(A_n)^{(n-1)}_{\sss H}$. 
\end{defn}


\begin{defn}\label{divisibletfsolvablegp} A torsion-free-solvable group $A$ is a \textbf{(uniquely) divisible torsion-free-solvable} group if, for each $n$, $A^{\sss (n)}_{\sss H}/A^{\sss (n+1)}_{\sss H}$ is a (uniquely) divisible $\mathbb{Z}[A/A^{\sss (n)}_{\sss H}]$-module.
\end{defn}

\begin{thm}\label{superharvey} For any group $G$ and any $n\geq 0$ there exist uniquely divisible $m$-torsion-free-solvable groups, $\wt G_m$, $0\leq m\leq n$, and a torsion-free-solvable tower, $\{\wt G_m, f_m:G\to\wt G_m, \pi_m:\wt G_m\to\wt G_{m-1}\}$, $0\leq m\leq n$ such that

\begin{enumerate}

\item $\ker f_m=G^{(m)}_H$.
\item If $A$ is finitely generated,
$B$ is finitely presented and $\phi:A\to B$ induces an isomorphism (respectively, monomorphism)
on $\hq1$ and induces an epimorphism $\phi_*: H_2(A;\mathbb{Q})\to H_2(B;\mathbb{Q})/<\Phi^{(n)}_H(B)>$ (that is, the cokernel of $\phi_*: H_2(A;\mathbb{Q})\to H_2(B;\mathbb{Q})$ is spanned by $B^{(n)}_H$-surfaces), then there is an
isomorphism (respectively, monomorphism) $\tilde\phi_n:\wt A_n\to\wt B_n$ such that the following
commutes
$$
\begin{array}{cccccc}
A     &\overset{f_n^A}{\lra}    &\wt A_n   &\overset{\pi_n^A}{\lra}
&\wt A_{n-1}\\
\phi\Big\downarrow         &&\tilde\phi_n\Big\downarrow
&&\Big\downarrow\tilde\phi_{n-1}\\
B   &\overset{f^B_n}{\lra}   &\wt B_n   &\overset{\pi^B_n}{\lra}
&\wt B_{n-1} \\
\end{array}
$$
\end{enumerate}
\end{thm}

The proof is identical to that in ~\cite{CH1}, with Proposition~\ref{2-connected} used in place of the weaker version used in ~\cite{CH1}.

\bibliographystyle{plain}
\bibliography{mybib4}

\end{document}

%% file: amsmath-defs.tex
\def\b{\beta}
\def\cd{\cdot}

\def\f{\frac}
\def\g{\gamma}
\def\G{\Gamma}

\def\lb\{{\left\{}
\def\la{\lambda}
\def\La{\Lambda}
\def\lla{\longleftarrow}
\def\lm{\limits}
\def\lra{\longrightarrow}
\def\dllra{\Longleftrightarrow}
\def\llra{\longleftrightarrow}
\def\n{\nabla}
\def\ngth{\negthickspace}
\def\ola{\overleftarrow}
\def\Om{\Omega}
\def\om{\omega}
\def\op{\oplus}
\def\oper{\operatorname}
\def\oplm{\operatornamewithlimits}
\def\ora{\overrightarrow}
\def\ov{\overline}
\def\ova{\overarrow}
\def\ox{\otimes}
\def\p{\partial}
\def\rb\}{\right\}}
\def\s{\sigma}
\def\sbq{\subseteq}
\def\spq{\supseteq}
\def\sqp{\sqsupset}
\def\supth{{\text{th}}}
\def\T{\Theta}
\def\th{\theta}
\def\tl{\tilde}
\def\thra{\twoheadrightarrow}
\def\un{\underline}
\def\ups{\upsilon}
\def\vp{\varphi}
\def\wh{\widehat}
\def\wt{\widetilde}
\def\x{\times}
\def\z{\zeta}
\def\({\left(}
\def\){\right)}
\def\[{\left[}
\def\]{\right]}
\def\<{\left<}
\def\>{\right>}

\def\tec{Teichm\"uller\ }
\def\sconr{\hbox{\medspace\vrule width 0.4pt height 4.7pt depth
0.4pt \vrule width 5pt height 0pt depth 0.4pt\medspace}}

\def\SA{\mathcal A}
\def\SB{\mathcal B}
\def\SC{\mathcal C}
\def\SD{\mathcal D}
\def\SE{\mathcal E}
\def\SF{\mathcal F}
\def\SG{\mathcal G}
\def\SH{\mathcal H}
\def\SI{\mathcal I}
\def\SJ{\mathcal J}
\def\SK{\mathcal K}
\def\SL{\mathcal L}
\def\SM{\mathcal M}
\def\SN{\mathcal N}
\def\SO{\mathcal O}
\def\SP{\mathcal P}
\def\SQ{\mathcal Q}
\def\SR{\mathcal R}
\def\SS{\mathcal S}
\def\ST{\mathcal T}
\def\SU{\mathcal U}
\def\SV{\mathcal V}
\def\SW{\mathcal W}
\def\SX{\mathcal X}
\def\SY{\mathcal Y}
\def\SZ{\mathcal Z}

\newcommand{\BA}{\ensuremath{\mathbf A}}
\newcommand{\BB}{\ensuremath{\mathbf B}}
\newcommand{\BC}{\ensuremath{\mathbf C}}
\newcommand{\BD}{\ensuremath{\mathbf D}}
\newcommand{\BE}{\ensuremath{\mathbf E}}
\newcommand{\BF}{\ensuremath{\mathbf F}}
\newcommand{\BG}{\ensuremath{\mathbf G}}
\newcommand{\BH}{\ensuremath{\mathbf H}}
\newcommand{\BI}{\ensuremath{\mathbf I}}
\newcommand{\BJ}{\ensuremath{\mathbf J}}
\newcommand{\BK}{\ensuremath{\mathbf K}}
\newcommand{\BL}{\ensuremath{\mathbf L}}
\newcommand{\BM}{\ensuremath{\mathbf M}}
\newcommand{\BN}{\ensuremath{\mathbf N}}
\newcommand{\BO}{\ensuremath{\mathbf O}}
\newcommand{\BP}{\ensuremath{\mathbf P}}
\newcommand{\BQ}{\ensuremath{\mathbf Q}}
\newcommand{\BR}{\ensuremath{\mathbf R}}
\newcommand{\BS}{\ensuremath{\mathbf S}}
\newcommand{\BT}{\ensuremath{\mathbf T}}
\newcommand{\BU}{\ensuremath{\mathbf U}}
\newcommand{\BV}{\ensuremath{\mathbf V}}
\newcommand{\BW}{\ensuremath{\mathbf W}}
\newcommand{\BX}{\ensuremath{\mathbf X}}
\newcommand{\BY}{\ensuremath{\mathbf Y}}
\newcommand{\BZ}{\ensuremath{\mathbf Z}}

\def\bba{{\mathbb A}}
\def\bbb{{\mathbb B}}
\def\bbc{{\mathbb C}}
\def\bbd{{\mathbb D}}
\def\bbe{{\mathbb E}}
\def\bbf{{\mathbb F}}
\def\bbg{{\mathbb G}}
\def\bbh{{\mathbb H}}
\def\bbi{{\mathbb I}}
\def\bbj{{\mathbb J}}
\def\bbk{{\mathbb K}}
\def\bbl{{\mathbb L}}
\def\bbm{{\mathbb M}}
\def\bbn{{\mathbb N}}
\def\bbo{{\mathbb O}}
\def\bbp{{\mathbb P}}
\def\bbq{{\mathbb Q}}
\def\bbr{{\mathbb R}}
\def\bbs{{\mathbb S}}
\def\bbt{{\mathbb T}}
\def\bbu{{\mathbb U}}
\def\bbv{{\mathbb V}}
\def\bbw{{\mathbb W}}
\def\bbx{{\mathbb X}}
\def\bby{{\mathbb Y}}
\def\bbz{{\mathbb Z}}